\documentclass{article}
\usepackage{array,amsmath,amsthm}
\numberwithin{equation}{section}
\usepackage{amssymb}
\usepackage{indentfirst}
\usepackage{geometry}
\begin{document}

\newtheorem{thm}{Theorem}[section]
\newtheorem{cor}[thm]{Corollary}
\newtheorem{prop}[thm]{Proposition}
\newtheorem{conj}[thm]{Conjecture}
\newtheorem{lem}[thm]{Lemma}
\newtheorem{Def}[thm]{Definition}
\newtheorem{rem}[thm]{Remark}
\newtheorem{prob}[thm]{Problem}
\newtheorem{ex}{Example}[section]

\newcommand{\be}{\begin{equation}}
\newcommand{\ee}{\end{equation}}
\newcommand{\ben}{\begin{enumerate}}
\newcommand{\een}{\end{enumerate}}
\newcommand{\beq}{\begin{eqnarray}}
\newcommand{\eeq}{\end{eqnarray}}
\newcommand{\beqn}{\begin{eqnarray*}}
\newcommand{\eeqn}{\end{eqnarray*}}
\newcommand{\bei}{\begin{itemize}}
\newcommand{\eei}{\end{itemize}}

\newcommand{\pa}{{\partial}}
\newcommand{\V}{{\rm V}}
\newcommand{\R}{{\bf R}}
\newcommand{\K}{{\rm K}}
\newcommand{\e}{{\epsilon}}
\newcommand{\tomega}{\tilde{\omega}}
\newcommand{\tOmega}{\tilde{Omega}}
\newcommand{\tR}{\tilde{R}}
\newcommand{\tB}{\tilde{B}}
\newcommand{\tGamma}{\tilde{\Gamma}}
\newcommand{\fa}{f_{\alpha}}
\newcommand{\fb}{f_{\beta}}
\newcommand{\faa}{f_{\alpha\alpha}}
\newcommand{\faaa}{f_{\alpha\alpha\alpha}}
\newcommand{\fab}{f_{\alpha\beta}}
\newcommand{\fabb}{f_{\alpha\beta\beta}}
\newcommand{\fbb}{f_{\beta\beta}}
\newcommand{\fbbb}{f_{\beta\beta\beta}}
\newcommand{\faab}{f_{\alpha\alpha\beta}}

\newcommand{\pxi}{ {\pa \over \pa x^i}}
\newcommand{\pxj}{ {\pa \over \pa x^j}}
\newcommand{\pxk}{ {\pa \over \pa x^k}}
\newcommand{\pyi}{ {\pa \over \pa y^i}}
\newcommand{\pyj}{ {\pa \over \pa y^j}}
\newcommand{\pyk}{ {\pa \over \pa y^k}}
\newcommand{\dxi}{{\delta \over \delta x^i}}
\newcommand{\dxj}{{\delta \over \delta x^j}}
\newcommand{\dxk}{{\delta \over \delta x^k}}

\newcommand{\px}{{\pa \over \pa x}}
\newcommand{\py}{{\pa \over \pa y}}
\newcommand{\pt}{{\pa \over \pa t}}
\newcommand{\ps}{{\pa \over \pa s}}
\newcommand{\pvi}{{\pa \over \pa v^i}}
\newcommand{\ty}{\tilde{y}}
\newcommand{\bGamma}{\bar{\Gamma}}

\title { Some inequalities and gradient estimates for harmonic functions on Finsler measure spaces  \footnote{The first author is supported by the National Natural Science Foundation of China (11871126, 12141101)}}
\author{ Xinyue Cheng$^{1}$ \& Yalu Feng$^{1}$ \footnote{Corresponding author}\\
$^{1}$ School of Mathematical Sciences, Chongqing Normal University, \\
Chongqing, 401331, P.R. China\\
E-mails: chengxy@cqnu.edu.cn \& fengyl2824@qq.com }
\date{}

\maketitle

\begin{abstract}
In this paper, we study functional and geometric inequalities on complete Finsler measure spaces under the condition that the weighted Ricci curvature ${\rm Ric}_\infty$ has a lower bound. We first obtain some local uniform Poincar\'{e} inequalities and Sobolev inequalities. Further, we give a mean value inequality for nonnegative subsolutions of elliptic equations. Finally, we obtain local and global Harnack inequalities, and then, establish a global gradient estimate for positive harmonic functions on forward complete non-compact Finsler measure spaces. Besides, as a by-product of the mean value inequality, we prove a Liouville type theorem. 
\end{abstract}

\maketitle

\section{Introduction}

The study of Harmonic functions on Riemannian manifolds is one of the important research topics in geometric analysis. It is well known that Yau's gradient estimate and Cheng-Yau's local gradient estimate for positive harmonic functions are important and essential results in Riemannian geometry \cite{CYau, Yau}. Latter, these results are generalized in different settings by many mathematicians. For example, Wang-Zhang obtained a local gradient estimate for positive $p(> 1)$-harmonic functions in a geodesic ball $B_R(o)$ on a Riemannian manifold $M$ \cite{XWang}. Munteanu-Wang generalized the above gradient estimate to the weighted Riemannian manifolds with nonnegative weighted Ricci curvature (that is, Bakry-\'{E}mery Ricci curvature) and derive a gradient estimate for positive $f$-harmonic functions \cite{MW1}. Zhang-Zhu generalized Yau's local gradient estimate to Alexandrov spaces \cite{HZhang}.

Finsler geometry is just Riemannian geometry without the quadratic restriction \cite{Chern}. Similar to Riemannian case, Finsler manifolds with Ricci curvature bounded below are always of some amazing properties. C. Xia generalized Cheng-Yau's result to Finsler manifolds and proved the local gradient estimate for harmonic functions on forward complete non-compact Finsler measure spaces under the condition that the weighted Ricci curvature ${\rm Ric}_{N}$ has a lower bound, ${\rm Ric}_{N}\geq-K$ for some $N\in [n, +\infty)$ and $K\geq 0$ \cite{CXia}. Recently, Q. Xia obtained the local and global gradient estimates for positive Finsler $p$-eigenfunctions on forward complete non-compact Finsler measure spaces under the condition that ${\rm Ric}_{N} \geq -K$ for some $N\in [n, \infty)$ and $K\geq 0$, and as applications, she obtained some Liouville and Harnack theorems and a global gradient estimate for positive $p$-harmonic functions \cite{QXia}.

On the other hand, compared with Riemannian case, it is natural to characterize functional and geometric inequalities on Finsler measure spaces under the condition about the weighted Ricci curvature ${\rm Ric}_\infty$. Actually, the role played by ${\rm Ric}_{N}$ and the role played by ${\rm Ric}_{\infty}$ in geometry and geometric analysis are usually quite different. Besides, in the studies of many problems, the  results under the condition about ${\rm Ric}_\infty$ can not be obtained from the corresponding results under the condition about ${\rm Ric }_N$ by letting $N \rightarrow\infty$.  Maybe one of the main reasons for these cases is that ${\rm Ric}_{\infty}$ can not fully control the volume properties of Finsler manifolds. For some details, see \cite{ChSh} and \cite{Zhu}.

In this paper, by $(M, F, m)$ we always denote a Finsler manifold $(M, F)$ equipped with a smooth measure $m$ which we call a Finsler measure space.
A Finsler measure space is not a metric space in usual sense because Finsler metric $F$ may be nonreversible, that is, $F(x, y)\neq F(x, -y)$ may happen.
This non-reversibility causes the asymmetry of the associated distance function. In order to overcome this defect, Ohta extended the concepts of uniform smoothness and the uniform convexity in Banach space theory into Finsler geometry and gave their geometric interpretation \cite{Ohta}.
The uniform smoothness and uniform convexity mean that there exist two uniform constants $0<\kappa^{*}\leq 1 \leq \kappa <\infty$ such that for $x\in M$, $V\in T_xM\setminus \{0\}$ and $W\in T_xM$, we have
\begin{equation}
\kappa^*F^2(x, W)\leq g_V(W, W)\leq \kappa F^2(x, W), \label{usk}
\end{equation}
where $g_V$ is the induced weighted Riemann metric on the tangent bundle of corresponding Finsler manifolds (see (\ref{weiRiem})). On the other hand, Rademacher defined the reversibility $\Lambda$ of $F$ by
\be
\Lambda:=\sup _{(x, y) \in TM \backslash\{0\}} \frac{F(x,y)}{F(x, -y)}.
\ee
Obviously, $\Lambda \in [1, \infty]$ and $\Lambda=1$ if and only if $F$ is reversible \cite{Ra}. If $F$ satisfies the uniform smoothness and uniform convexity, then $\Lambda$ is finite with
$$
1 \leq \Lambda \leq \min \left\{\sqrt{\kappa}, \sqrt{1 / \kappa^*}\right\}.
$$
$F$ is Riemannian if and only if $\kappa=1$ if and only if $\kappa^*=1$ \cite{ChernShen, Ohta, Ra}.

Let $x_{0}\in M$. The forward and backward geodesic balls of radius $R$ with center at $x_{0}$ are respectively defined by
$$
B_R^{+}(x_0):=\{x \in M \mid d(x_{0}, x)<R\},\qquad B_R^{-}(x_0):=\{x \in M \mid d(x, x_{0})<R\}.
$$
In the following, we will always denote $B_{R}:=B^{+}_{R}(x_0)$ for some $x_{0}\in M$ for simplicity. For more details, see Section \ref{Introd}.

Our initial motivation and objective is to prove the global gradient estimate for harmonic functions on forward complete non-compact Finsler measure spaces under the condition that the weighted Ricci curvature ${\rm Ric}_{\infty}$ has a lower bound. For this aim, we will first derive some volume comparison theorems of Bishop-Gromov type based on Laplacian comparison theorems given by Z. Shen and S. Yin in \cite{Shen2} and \cite{Yin} respectively (also see \cite{ChSh}) in Section \ref{VolCom}. Then, as the applications of volume comparison theorems and for the need to run Moser's iteration, we will give a local uniform Poincar\'{e} inequality and a local uniform Sobolev inequality  by following the similar argument of Lemma 3.2 in  Munteanu-Wang \cite{MW1}  under the condition that ${\rm Ric}_{\infty}$ has a lower bound in Section \ref{PSIne}. Further, by combining these obtained inequalities, we will prove mean value inequality for nonnegative subsolutions of elliptic equations and Harnack inequality of harmonic functions in Section \ref{proof}. Finally, based on the mean value inequality  and the Harnack inequality, we will prove the gradient estimate that we want to get still in Section  \ref{proof}.

In this paper, we always assume that the dimension $n$ of the manifolds satisfies that $n\geq2$. The first main result in this paper is the following theorem which gives a mean value inequality for nonnegative subsolutions of elliptic equations.
\begin{thm}\label{meanineq}
Let $(M, F, m)$ be an $n$-dimensional forward complete Finsler measure space with finite reversibility $\Lambda$. Assume that ${\rm Ric}_{\infty}\geq K$ and $|\tau|\leq k$ for some $K\in \mathbb{R}$ and $k>0$. Suppose that $u$ is a nonnegative function defined on $B_R$ satisfying
$$
\Delta u \geq -fu
$$
in the weak sense, where $f\in L^{\infty}(B_R)$ is nonnegative. Here, $\tau$ denotes the distortion of $F$. Then for any $p \in(0,2]$ and $\delta \in(0,1)$, there are constants $\nu>2$ and $C=C(n, k, \nu, p, \Lambda)>0$ depending on $n, k, \nu, p$ and  $\Lambda$, such that
\begin{equation}
\sup _{B_{\delta R}} u^p \leq e^{C(1+\sqrt{|K|} R)}(1+\mathcal{A}R^2)^{\frac{\nu}{2}} (1-\delta)^{-\nu} m\left(B_R\right)^{-1}  \int_{B_R} u^{p} dm, \label{meanineq-1}
\end{equation}
where $\mathcal{A}:=\sup\limits_{B_R}|f|$. In particular, if ${\rm Ric}_{\infty}\geq 0$, we have
\begin{equation}
\sup _{B_{\delta R}} u^p \leq C (1+\mathcal{A}R^2)^{\frac{\nu}{2}} (1-\delta)^{-\nu} m\left(B_R\right)^{-1} \int_{B_R} u^{p} dm. \label{meanineq-2}
\end{equation}
\end{thm}

As an application of Theorem \ref{meanineq}, we obtain the following Liouville property of harmonic function.

\begin{cor}\label{Liouville}
Let $(M, F, m)$ be an $n$-dimensional forward complete and non-compact Finsler measure space equipped with a uniformly convex and uniformly smooth Finsler metric $F$. Assume that ${\rm Ric}_{\infty}\geq 0$ and $|\tau|\leq k$ for some $k>0$. Then any positive harmonic function $u$ satisfying
\begin{equation}
\lim\limits_{r\rightarrow\infty}\frac{|u(x)|}{r}=0    \label{L-1}
\end{equation}
must be a constant, where $r= d(x_{0}, x)$ denotes the distance function from some point $x_{0} \in M$.
\end{cor}

In fact, Corollary \ref{Liouville} can be regarded as an analogue under the condition that ${\rm Ric}_{\infty}\geq 0$ of Corollary 1.3 in \cite{CXia}, which was obtained under the condition that ${\rm Ric}_{N} \geq 0$.

It should be pointed out that the mean value inequality in Theorem \ref{meanineq} is very important in the following study of Harnack inequality and gradient estimates for harmonic functions on Finsler measure spaces.  From Theorem \ref{meanineq}, we can obtain the following Harnack inequality for positive harmonic functions.

\begin{thm}{\label{mean2}}
Let $(M, F, m)$ be an $n$-dimensional forward complete Finsler measure space with finite reversibility $\Lambda$. Assume that ${\rm Ric}_{\infty} \geq K$ and $|\tau|\leq k$ for some $K\in \mathbb{R}$  and $k>0$. If $u$ is a positive harmonic function on $B_R$, then,  for any $\delta \in (0,1)$,  there exist positive constant $C=C\left(n, \delta, k, \Lambda\right)$ depending on $n, \delta, k$ and  $\Lambda$, such that
\begin{equation}
\sup\limits_{B_{\delta R}}u \leq e^{C(1+\sqrt{|K|} R)}\inf\limits_{B_{\delta R}}u. \label{mean-1}
\end{equation}
\end{thm}

In \cite{CXia}, C. Xia obtained the Harnack inequality as an application of the local gradient estimate for positive harmonic functions on a forward complete non-compact Finsler measure space under the condition that ${\rm Ric}_{N}\geq -K$ for some $N\in [n, +\infty)$ and $K\geq 0$. Further, under the same conditions, Q. Xia obtained the Harnack inequality as an application of local gradient estimate for positive $p$($>1$)-eigenfunctions or $p$-harmonic functions \cite{QXia}.

It is worth to note that C. Xia proved the local gradient estimate for positive harmonic functions on a forward complete non-compact Finsler measure space by local uniform Sobolev inequality and Moser's iteration under the condition that ${\rm Ric}_{N}\geq -K$ for $N\in [n, +\infty)$ and $K \geq 0$ \cite{CXia}. Also, Q. Xia proved the local gradient estimate for positive $p$-eigenfunctions on a forward complete non-compact Finsler measure space by same arguments under the same conditions \cite{QXia}. However, under the condition that Ric$_\infty\geq -K$ for some $K\geq0$, it is difficult to obtain the  gradient estimate for positive harmonic functions on a forward complete non-compact Finsler measure space by using similar methods, even if one adds the condition about distortion $\tau$. Therefore, we will follow the standard arguments in Munteanu-Wang \cite{MW2} to prove the following gradient estimate for positive harmonic functions, which means that we will get our gradient estimate via the mean value inequality and the Harnack inequality for positive harmonic functions.

\begin{thm}\label{ggeh}
Let $(M, F, m)$ be an $n$-dimensional forward complete non-compact Finsler measure space equipped with a uniformly convex and uniformly smooth Finsler metric $F$. Assume that ${\rm Ric}_{\infty}\geq -K$ and $|\tau|\leq k$ for some $K \geq 0$ and $k>0$. Let $u$ be a positive harmonic
function on $M$, that is, $\Delta u=0$
in a weak sense on $M$. Then there exists a positive constant $C=C\left(n, k, \kappa, \kappa^*\right)$ depending only on $n, k$, the uniform constants $\kappa$ and $\kappa^*$, such that
\begin{equation}
\max\limits_{x\in M}\{F(x, \nabla \log u(x)), F(x, \nabla(-\log u(x)))\} \leq e^{C(1+\sqrt{K})}(1+K)^{n+4k}. \label{gradient-1}
\end{equation}
\end{thm}

As a special case of Theorem \ref{ggeh}, the following result can be regarded as a direct generalization of Munteanu-Wang's result in \cite{MW2} under the condition that $|\tau|\leq k$ on Finsler measure spaces.

\begin{cor}\label{geh}
Let $(M, F, m)$ be an $n$-dimensional forward complete non-compact Finsler measure space equipped with a uniformly convex and uniformly smooth Finsler metric $F$. Assume that ${\rm Ric}_{\infty}\geq-(n-1)$ and $|\tau|\leq k$ for some $k>0$. Let $u$ be a positive harmonic function on $M$, i.e., $\Delta u=0$
in a weak sense on $M$. Then there exists a positive constant $C=C\left(n, k, \kappa, \kappa^*\right)$ depending on $n$, $k$, the uniform constants $\kappa$ and $\kappa^*$, such that
\begin{equation}\label{gradient-2}
\max\limits_{x\in M}\{F(x, \nabla \log u(x)), F(x, \nabla(-\log u(x)))\} \leq C.
\end{equation}
\end{cor}

\section{Preliminaries}\label{Introd}
In this section, we briefly review some necessary definitions, notations and  fundamental results in Finsler geometry. For more details, we refer to \cite{BaoChern, ChernShen, OHTA, Shen1}.

\subsection{Finsler measure spaces}
Let $M$ be an $n$-dimensional smooth manifold. For a point $x \in M$, denote by $T_x M$ the tangent space of $M$ at $x$. The tangent bundle $T M$ of $M$ is the union of tangent spaces with a natural differential structure,
$$
T M=\bigcup_{x \in M} T_x M .
$$
Denote the elements in $T M$ by $(x, y)$ with $y \in T_x M$. Let $T M_0:=T M \backslash\{0\}$ and $\pi: T M \backslash\{0\} \rightarrow M$ be the natural projective map. The pull-back $\pi^* T M$ is an n-dimensional vector bundle on $T M_0$. A Finsler metric on manifold $M$ is a function $F: T M \longrightarrow[0, \infty)$ on the tangent bundle satisfying the following properties:
\ben
\item[{\rm (1)}]  $F$ is $C^{\infty}$ on $TM\backslash\{0\}$;
\item[{\rm (2)}]  $F(x,\lambda y)=\lambda F(x,y)$ for any $(x,y)\in TM$ and all $\lambda >0$;
\item[{\rm (3)}]  $F$ is strongly convex, that is, the matrix $\left(g_{ij}(x,y)\right)=\left(\frac{1}{2}(F^{2})_{y^{i}y^{j}}\right)$ is positive definite for any nonzero $y\in T_{x}M$.
\een
Such a pair $(M,F)$ is called a Finsler manifold and $g:=g_{ij}(x,y)dx^{i}\otimes dx^{j}$ is called the fundamental tensor of $F$.

We define the reverse metric $\overleftarrow{F} $ of $F$ by $\overleftarrow{F}(x, y):=F(x,-y)$ for all $(x, y) \in T M$. It is easy to see that $\overleftarrow{F}$ is also a Finsler metric on $M$. A Finsler metric $F$ on $M$ is said to be reversible if $\overleftarrow{F}(x, y)=F(x, y)$ for all $(x, y) \in T M$. Otherwise, we say $F$ is irreversible. For a non-vanishing vector field $V$ on $M$, one introduces the weighted Riemannian metric $g_V$ on $M$ given by
\be
g_V(y, w)=g_{ij}(x, V_x)y^i w^j  \label{weiRiem}
\ee
for $y,\, w\in T_{x}M$. In particular, $g_V(V,V)=F^2(V,V)$.

Let $(M,F)$ be a Finsler manifold of dimension $n$. The pull-back $\pi ^{*}TM$ admits a unique linear connection, which is called the Chern connection. The Chern connection $D$ is determined by the following equations
\beq
&& D^{V}_{X}Y-D^{V}_{Y}X=[X,Y], \label{chern1}\\
&& Zg_{V}(X,Y)=g_{V}(D^{V}_{Z}X,Y)+g_{V}(X,D^{V}_{Z}Y)+ 2C_{V}(D^{V}_{Z}V,X,Y) \label{chern2}
\eeq
for $V\in TM\setminus \{0\}$  and $X, Y, Z \in TM$, where
$$
C_{V}(X,Y,Z):=C_{ijk}(x,V)X^{i}Y^{j}Z^{k}=\frac{1}{4}\frac{\pa ^{3}F^{2}(x,V)}{\pa V^{i}\pa V^{j}\pa V^{k}}X^{i}Y^{j}Z^{k}
$$
is the Cartan tensor of $F$ and $D^{V}_{X}Y$ is the covariant derivative with respect to the reference vector $V$.

For $x_1, x_2 \in M$, the distance from $x_1$ to $x_2$ is defined by
$$
d\left(x_1, x_2\right):=\inf _\gamma \int_0^1 F(\dot{\gamma}(t)) d t,
$$
where the infimum is taken over all $C^1$ curves $\gamma:[0,1] \rightarrow M$ such that $\gamma(0)=$ $x_1$ and $\gamma(1)=x_2$. Note that $d \left(x_1, x_2\right) \neq d \left(x_2, x_1\right)$ unless $F$ is reversible.
A $C^{\infty}$-curve $\gamma:[0,1] \rightarrow M$ is called a geodesic  if $F(\gamma, \dot{\gamma})$ is constant and it is locally minimizing.

The exponential map $\exp _x: T_x M \rightarrow M$ is defined by $\exp _x(v)=\gamma(1)$ for $v \in T_x M$ if there is a geodesic $\gamma:[0,1] \rightarrow M$ with $\gamma(0)=x$ and $\dot{\gamma}(0)=v$. A Finsler manifold $(M, F)$ is said to be forward complete (resp. backward complete) if each geodesic defined on $[0, \ell)$ (resp. $(-\ell, 0])$ can be extended to a geodesic defined on $[0, \infty)$ (resp. $(-\infty, 0])$. We say $(M, F)$ is complete if it is both forward complete and backward complete. By Hopf-Rinow theorem on forward complete Finsler manifolds, any two points in $M$ can be connected by a minimal forward geodesic and the forward closed balls $\overline{B_R^{+}(p)}$ are compact. For a point $p \in M$ and a unit vector $v \in T_p M$, let $\rho(v):=\sup \left\{t>0 \mid \text{the geodesic} \ \exp _p(tv) \ \text{is minimal} \right\}$. If $\rho(v)<\infty$, we call $\exp _p\left(\rho(v) v\right)$ a cut point of $p$. The set of all the cut points of $p$ is said to be the cut locus of $p$, denoted by $Cut(p)$. The cut locus of $p$ always has null measure and $d_{p}:= d (p, \cdot)$ is $C^1$ outside the cut locus of $p$ (see \cite{BaoChern, Shen1}).

Let $(M, F, m)$ be an $n$-dimensional Finsler manifold with a smooth measure $m$. Write the volume form $dm$ of  $m$ as $d m = \sigma(x) dx^{1} dx^{2} \cdots d x^{n}$. Define
\be\label{Dis}
\tau (x, y):=\ln \frac{\sqrt{{\rm det}\left(g_{i j}(x, y)\right)}}{\sigma(x)}.
\ee
We call $\tau$ the distortion of $F$. It is natural to study the rate of change of the distortion along geodesics. For a vector $y \in T_{x} M \backslash\{0\}$, let $\sigma=\sigma(t)$ be the geodesic with $\sigma(0)=x$ and $\dot{\sigma}(0)=y.$  Set
\be
{\bf S}(x, y):= \frac{d}{d t}\left[\tau(\sigma(t), \dot{\sigma}(t))\right]|_{t=0}.
\ee
$\mathbf{S}$ is called the S-curvature of $F$ \cite{ChernShen, shen}.

\subsection{Gradient and Finsler Laplacian}

Given a Finsler structure $F$ on $M$,  there is a Finsler co-metric $F^{*}$ on $M$ which is non-negative function on the cotangent bundle $T^{*}M$ given by
\be
F^{*}(x, \xi):=\sup\limits_{y\in T_{x}M\setminus \{0\}} \frac{\xi (y)}{F(x,y)}, \ \ \forall \xi \in T^{*}_{x}M. \label{co-Finsler}
\ee
We call $F^{*}$ the dual Finsler metric of $F$.  For any vector $y\in T_{x}M\setminus \{0\}$, $x\in M$, the covector $\xi =g_{y}(y, \cdot)\in T^{*}_{x}M$ satisfies
\be
F(x,y)=F^{*}(x, \xi)=\frac{\xi (y)}{F(x,y)}. \label{shenF311}
\ee
Conversely, for any covector $\xi \in T_{x}^{*}M\setminus \{0\}$, there exists a unique vector $y\in T_{x}M\setminus \{0\}$ such that $\xi =g_{y}(y, \cdot)\in T^{*}_{x}M$ (Lemma 3.1.1, \cite{Shen1}). Naturally,  we define a map ${\cal L}: TM \rightarrow T^{*}M$ by
$$
{\cal L}(y):=\left\{
\begin{array}{ll}
g_{y}(y, \cdot), & y\neq 0, \\
0, & y=0.
\end{array} \right.
$$
It follows from (\ref{shenF311}) that
$$
F(x,y)=F^{*}(x, {\cal L}(y)).
$$
Thus ${\cal L}$ is a norm-preserving transformation. We call ${\cal L}$ the Legendre transformation on Finsler manifold $(M, F)$.

Let
$$
g^{*kl}(x,\xi):=\frac{1}{2}\left[F^{*2}\right]_{\xi _{k}\xi_{l}}(x,\xi).
$$
For any $\xi ={\cal L}(y)$, we have
\be
g^{*kl}(x,\xi)=g^{kl}(x,y), \label{Fdual}
\ee
where $\left(g^{kl}(x,y)\right)= \left(g_{kl}(x,y)\right)^{-1}$. If $F$ is uniformly smooth and convex with (\ref{usk}), then  $\left(g^{ij}\right)$ is uniformly elliptic in the sense that there exists two constants $\tilde{\kappa}=(\kappa^*)^{-1}$, $\tilde{\kappa}^*=\kappa^{-1}$ such that for $x \in M, \ \xi \in T^*_x M \backslash\{0\}$ and $\eta \in T_x^* M$, we have
\be
\tilde{\kappa}^* F^{* 2}(x, \eta) \leq g^{*i j}(x, \xi) \eta_i \eta_j \leq \tilde{\kappa} F^{* 2}(x, \eta). \label{unisc}
\ee

Given a smooth function $u$ on $M$, the differential $d u_x$ at any point $x \in M$,
$$
d u_x=\frac{\partial u}{\partial x^i}(x) d x^i
$$
is a linear function on $T_x M$. We define the gradient vector $\nabla u(x)$ of $u$ at $x \in M$ by $\nabla u(x):=\mathcal{L}^{-1}(d u(x)) \in T_x M$. In a local coordinate system, we can express $\nabla u$ as
\be \label{nabna}
\nabla u(x)= \begin{cases}g^{* i j}(x, d u) \frac{\partial u}{\partial x^i} \frac{\partial}{\partial x^j}, & x \in M_u, \\ 0, & x \in M \backslash M_u,\end{cases}
\ee
where $M_{u}:=\{x \in M \mid d u(x) \neq 0\}$ \cite{Shen1}. In general, $\nabla u$ is only continuous on $M$, but smooth on $M_{u}$.

The Hessian of $u$ is defined by using Chern connection as
$$
\nabla^2 u(X, Y)=g_{\nabla u}\left(D_X^{\nabla u} \nabla u, Y\right).
$$
One can show that $\nabla^2 u(X, Y)$ is symmetric, see \cite{Ohta2, WuXin}.

Let $(M, F, m)$ be an $n$-dimensional Finsler manifold with a smooth measure $m$. We may decompose the volume form $d m$ of $m$ as $d m=\mathrm{e}^{\Phi} d x^1 d x^2 \cdots d x^n$. Then the divergence of a differentiable vector field $V$ on $M$ is defined by
$$
\operatorname{div}_m V:=\frac{\partial V^i}{\partial x^i}+V^i \frac{\partial \Phi}{\partial x^i}, \quad V=V^i \frac{\partial}{\partial x^i} .
$$
One can also define $\operatorname{div}_m V$ in the weak form by following divergence formula
$$
\int_M \phi \operatorname{div}_m V d m=-\int_M d \phi(V) d m
$$
for all $\phi \in \mathcal{C}_0^{\infty}(M)$. Now we define the Finsler Laplacian $\Delta u$ by
\be
\Delta u:=\operatorname{div}_m(\nabla u). \label{Lap}
\ee
From (\ref{Lap}), Finsler Laplacian is a nonlinear elliptic differential operator of the second order.

Let $W^{1, p}(M)(p>1)$ be the space of functions $u \in L^p(M)$ with $\int_M[F(\nabla u)]^p d m<\infty$ and $W_0^{1, p}(M)$ be the closure of $\mathcal{C}_0^{\infty}(M)$ under the (absolutely homogeneous) norm
\be
\|u\|_{W^{1, p}(M)}:=\|u\|_{L^p(M)}+\frac{1}{2}\|F(\nabla u)\|_{L^p(M)}+\frac{1}{2}\|\overleftarrow{F}(\overleftarrow{\nabla} u)\|_{L^p(M)},
\ee
where $\overleftarrow{\nabla} u$ is the gradient of $u$ with respect to the reverse metric $\overleftarrow{F}$. In fact $\overleftarrow{F}(\overleftarrow{\nabla} u)=F(\nabla(-u))$.

Note that $\nabla u$ is weakly differentiable, the Finsler Laplacian should be understood in a weak sense, that is, for $u \in W^{1,2}(M)$, $\Delta u$ is defined by
\be
\int_M \phi \Delta u d m:=-\int_M d \phi(\nabla u) dm
\ee
for $\phi \in \mathcal{C}_0^{\infty}(M)$ \cite{Shen1}.

Given a weakly differentiable function $u$ and a vector field $V$ which does not vanish on $M_u$, the weighted Laplacian of $u$ on the weighted Riemannian manifold $\left(M, g_V, m\right)$ is defined  by
$$
\Delta^{V} u:= {\rm div}_{m}\left(\nabla^V u\right),
$$
where
$$
\nabla^V u:= \begin{cases}g^{ij}(x, V) \frac{\partial u}{\partial x^i} \frac{\partial}{\partial x^j} & \text { for } x \in M_u, \\ 0 & \text { for } x \notin M_u .\end{cases}
$$
Similarly, the weighted Laplacian can be viewed in a weak sense. We note that $\nabla^{\nabla u}u=\nabla u$ and $\Delta^{\nabla u} u=$ $\Delta u$. Moreover, it is easy to see that $\Delta u= {\rm tr}_{\nabla u} \nabla^2 u-{\bf S}(\nabla u)$ on $M_u$ \cite{OHTA, WuXin}.

\subsection{Weighted Ricci curvature}

Let $(M, F, m)$ be an $n$-dimensional  Finsler measure space. Given a non-vanishing vector field $V$ on $M$,  the Riemannian curvature  $R^V$ is defined by
$$
R^V(X, Y) Z=D_X^V D_Y^V Z-D_Y^V D_X^V Z-D_{[X, Y]}^V Z
$$
for any vector fields $X$, $Y$, $Z$ on $M$. For two linearly independent vectors $V, W \in T_x M \backslash\{0\}$, the flag curvature is defined by
$$
\mathcal{K}^V(V, W)=\frac{g_V\left(R^V(V, W) W, V\right)}{g_V(V, V) g_V(W, W)-g_V(V, W)^2}.
$$
Then the Ricci curvature is defined as
$$
\operatorname{Ric}(V):=F(x, V)^{2} \sum_{i=1}^{n-1} \mathcal{K}^V\left(V, e_i\right),
$$
where $e_1, \ldots, e_{n-1}, \frac{V}{F(V)}$ form an orthonormal basis of $T_x M$ with respect to $g_V$.

Let $Y$ be a $C^{\infty}$ geodesic field on an open subset $U \subset M$ and $\hat{g}=g_{Y}.$  Let
\be
d m:=e^{- \psi} {\rm Vol}_{\hat{g}}, \ \ \ {\rm Vol}_{\hat{g}}= \sqrt{{det}\left(g_{i j}\left(x, Y_{x}\right)\right)}dx^{1} \cdots dx^{n}. \label{voldecom}
\ee
It is easy to see that $\psi$ is given by
$$
\psi (x)= \ln \frac{\sqrt{\operatorname{det}\left(g_{i j}\left(x, Y_{x}\right)\right)}}{\sigma(x)}=\tau\left(x, Y_{x}\right),
$$
which is just the distortion along $Y_{x}$ at $x\in M$ \cite{ChernShen, Shen1}. Let $y := Y_{x}\in T_{x}M$ (that is, $Y$ is a geodesic extension of $y\in T_{x}M$). Then, by the definitions of the S-curvature, we have
\beqn
&&  {\bf S}(x, y)= Y[\tau(x, Y)]|_{x} = d \psi (y),  \\
&&  \dot{\bf S}(x, y)= Y[{\bf S}(x, Y)]|_{x} =y[Y(\psi)],
\eeqn
where $\dot{\bf S}(x, y):={\bf S}_{|m}(x, y)y^{m}$ and ``$|$" denotes the horizontal covariant derivative with respect to the Chern connection  \cite{shen, Shen1}. Further, the weighted Ricci curvatures are defined as follows \cite{ChSh,OHTA, Shen2}
\beq
{\rm Ric}_{N}(y)&=& {\rm Ric}(y)+ \dot{\bf S}(x, y) -\frac{{\bf S}(x, y)^{2}}{N-n},   \label{weRicci3}\\
{\rm Ric}_{\infty}(y)&=& {\rm Ric}(y)+ \dot{\bf S}(x, y). \label{weRicciinf}
\eeq
We say that Ric$_N\geq K$ for some $K\in \mathbb{R}$ if Ric$_N(v)\geq KF^2(v)$ for all $v\in TM$, where $N\in \mathbb{R}\setminus \{n\}$ or $N= \infty$.

The following Bochner-Weitzenb\"{o}ck type formula established by Ohta-Sturm \cite{Ohta2} is very important to derive the gradient estimates in this paper.

\begin{thm}{\rm (\cite{Ohta2,OHTA})}\label{boch} Given $u \in W_{\mathrm{loc}}^{2,2}(M) \bigcap C^1(M)$ with $\Delta u \in W_{\mathrm{loc}}^{1,2}(M)$, we have
\beq
&& -\int_M d \phi\left(\nabla^{\nabla u}\left(\frac{F^2(x, \nabla u)}{2}\right)\right) d m  \nonumber \\
&& \ \ =\int_M \phi \left\{d \left(\Delta u\right)(\nabla u)+{\rm Ric}_{\infty}(\nabla u)+\left\|\nabla^2 u\right\|_{H S(\nabla u)}^2\right\} dm \label{BWforinf}
\eeq
as well as
\be
-\int_M d \phi\left(\nabla^{\nabla u}\left(\frac{F^2(x, \nabla u)}{2}\right)\right) dm \geq \int_M \phi \left\{d \left(\Delta u\right)(\nabla u)+R i c_N(\nabla u)+\frac{\left(\Delta u\right)^2}{N}\right\} dm
\ee
for any $N \in[n, \infty]$ and all nonnegative functions $\phi \in W_0^{1,2}(M) \bigcap L^{\infty}(M)$. Here $\left\|\nabla^2 u\right\|_{H S(\nabla u)}^2$ denotes the Hilbert-Schmidt norm with respect to $g_{\nabla u}$.
\end{thm}

\section{Volume comparison theorems}\label{volcomth}\label{VolCom}

Let $(M, F, m)$ be an $n$-dimensional Finsler manifold with a smooth measure $m$ and $x \in M$. Let $\mathcal{D}_x:=M \backslash(\{x\} \cup C u t(x))$ be the cut-domain on $M$. For any $z \in \mathcal{D}_x$, we can choose the geodesic polar coordinates $(r, \theta)$ centered at $x$ for $z$ such that $r(z)=F(v)$ and $\theta^{\alpha}(z)=\theta^{\alpha}\left(\frac{v}{F(v)}\right)$, where $r(z)=d(x, z)$  and $v=\exp _x^{-1}(z) \in T_x M \backslash\{0\}$. It is well known that the distance function $r$ starting from $x \in M$ is smooth on $\mathcal{D}_x$ and $F(\nabla r)=1$ \cite{BaoChern,Shen1}.  A basic fact is that the distance function $r=d(x, \cdot)$ satisfies the following \cite{Shen1, Shen2}
$$
\nabla r |_{z}= \frac{\pa}{\pa r}|_{z}.
$$
By Gauss's lemma, the unit radial coordinate vector $\frac{\partial }{{\partial r}}$ and the coordinate vectors $\frac{\partial }{{\partial {\theta ^\alpha }}}$ for $1\leq \alpha \leq n-1$ are mutually vertical with respect to $g_{\nabla r}$ (Lemma 6.1.1 in \cite{BaoChern}).  Therefore, we can simply write the volume form at $z=\exp _{x}(r\xi)$  with $v=r \xi$ as $\left.dm\right|_{\exp _x(r \xi)}=\sigma(x, r, \theta) dr d\theta$, where $\xi \in I_{x}:=\left\{\xi \in T_x M \mid F(x, \xi)=1\right\}$.
Then, for forward geodesic ball $B_{R}=B_R^{+}(x)$ of radius $R$ at the center $x \in M$, the volume of $B_R$ is
$$
m(B_R)=\int_{B_R} d m=\int_{B_R \cap \mathcal{D}_x} d m=\int_0^{R} dr \int_{\mathcal{D}_x(r)} \sigma(x, r, \theta) d\theta,
$$
where $\mathcal{D}_x(r)=\left\{\xi \in I_x \mid r \xi \in \exp _x^{-1}\left(\mathcal{D}_x \cap B_R^{+}(x)\right)\right\}$. Obviously, for any $0<s<t<R, \ \mathcal{D}_x(t) \subseteq \mathcal{D}_x(s)$. Besides, by the definition \ref{Lap} of Laplacian, we have \cite{Shen1, Shen2, WuXin}
\be
\Delta r=\frac{\partial}{\partial r}\log \sigma (x, r, \theta). \label{Lapdis}
\ee

Define
$$
{s_c}(t): = \left\{ {\begin{array}{*{20}{c}}
\begin{array}{l}
\frac{1}{{\sqrt c }}\sin (\sqrt c t)\\
t\\
\frac{1}{{\sqrt { - c} }}\sinh (\sqrt { - c} t)
\end{array}&\begin{array}{l}
c > 0,\\
c = 0,\\
c < 0.
\end{array}
\end{array}} \right.
$$
Obviously, ${s_c}(t)$ is the solution of the differential equation $f'' + cf = 0$ satisfying $f(0) = 0$ and $f'(0) = 1$. Further, we define
$$
c{t_c}(t): = \frac{{{{s'}_c}(t)}}{{{s_c}(t)}} = \left\{ {\begin{array}{*{20}{c}}
\begin{array}{l}
\sqrt c \cot (\sqrt c t)\\
\frac{1}{t}\\
\sqrt { - c} \coth (\sqrt { - c} t)
\end{array}&\begin{array}{l}
c > 0,\\
c = 0,\\
c < 0.
\end{array}
\end{array}} \right.
$$

For any $p\in M$, let $r (x):=d (p,x)$ be the the distance function starting from $p$. We denote the forward geodesic sphere of radius $R$ at the center $p \in M$ by $S_{R}(p)$.

\begin{lem}{\rm(Laplacian comparison, \cite{ChSh,Shen2,Yin})}\label{Lapcompa}
Let $(M, F, m)$ be an $n$-dimensional Finsler manifold with weighted Ricci curvature satisfying $\operatorname{Ric}_{\infty} \geq K$. Then the following bound on $\Delta r$ holds.
\ben
\item[(1)] If distance function $r(x)$ is smooth and $r(x)>r_{0}$, then
$$
\Delta r \leq\left.\frac{d}{d t}[\ln \chi(t)]\right|_{t=r(x)},
$$
where $\chi(t)=e^{m_{0}\left(t-r_{0}\right)-\frac{1}{2} K\left(t-r_{0}\right)^{2}}\left(r_{0}<t<\infty\right)$ and $m_{0}:=\sup _{x \in r^{-1}\left(r_{0}\right)} \Delta r(x)$ characterizes the mean curvature of the geodesic sphere $S_{r_{0}}(p)$ {\rm \cite{Shen1}}.

\item[(2)] If $\mathbf{S}\geq-\alpha$, $\alpha>0$, then the following holds on $\mathcal{D}_{p} \cap B_{r_{o}}(p)$,
$$
\Delta r \leq\left.\frac{d}{d t}[\ln \chi(t)]\right|_{t=r(x)},
$$
where $\chi(t)=\left[s_{K/(n-1)}(t)\right]^{n-1} e^{\alpha t} \left(0<t<r_{o}\right)$. Here $r_{o}=+\infty$ when $K\leq 0$ and $r_{o}=\frac{\pi}{2} \sqrt{\frac{n-1}{K}}$ when $K>0$.

\item[(3)] If the distortion $\tau$ satisfies that $|\tau| \leq k$ for some $k>0$, then the following holds on $\mathcal{D}_{p} \cap B_{r_{o}}(p)$
$$
\Delta r \leq\left.\frac{d}{d t}[\ln \chi(t)]\right|_{t=r(x)},
$$
where $\chi(t):=\left[s_{K/(n-1)}(t)\right]^{n+4 k-1} \left(0<t<r_{o}\right)$. Here $r_{o}:=+\infty$ when $K\leq 0$ and $r_{o}=\frac{\pi}{4} \sqrt{\frac{n-1}{K}}$ when $K>0$.
\een
\end{lem}

From Lemma \ref{Lapcompa}, we can obtain the following volume comparison theorem.

\begin{prop}\label{volcom1} Let $(M, F, m)$ be an n-dimensional Finsler manifold with a smooth volume form dm. Assume that ${\rm Ric}_{\infty}\geq -K$, $K\geq 0$. Then, along any minimizing geodesic starting from the center $x$ of $B_{R}(x)$, the following conclusions hold.
\ben
\item[(1)] If $\mathbf{S}\geq -\alpha$, $\alpha>0$, then, for any $0<r_1<r_2<R$, we have
\be
\frac{\sigma\left(x, r_2, \theta\right)}{\sigma\left(x, r_1, \theta\right)} \leq\left(\frac{r_2}{r_1}\right)^{n-1} e^{r_2 (\alpha+\sqrt{(n-1)K})}. \label{ssigma}
\ee
Further, we have
\be
\frac{m\left(B_{r_2}(x)\right)}{m\left(B_{r_1}(x)\right)} \leq\left(\frac{r_2}{r_1}\right)^n e^{r_2(\alpha+\sqrt{(n-1)K})}. \label{volBallS}
\ee
\item[(2)] If $|\tau|\leq k$, $k>0$, then, for any $0<r_1<r_2<R$, we have
\be
\frac{\sigma\left(x, r_2, \theta\right)}{\sigma\left(x, r_1, \theta\right)} \leq\left(\frac{r_2}{r_1}\right)^{n+4k-1} e^{r_2(n+4k-1)\sqrt{K/(n-1)}}.
\ee
Further,  we have
\be
\frac{m\left(B_{r_2}(x)\right)}{m\left(B_{r_1}(x)\right)} \leq\left(\frac{r_2}{r_1}\right)^{n+4k} e^{r_2(n+4k)\sqrt{K/(n-1)}}.
\ee
\een
\end{prop}
{\it Proof.}  Firstly, let $\eta:[0, r] \rightarrow M$ be the minimizing geodesic from $\eta(0)=x$ to $\eta(r)=z$, where $r=d(x, z)$. By using the geodesic polar coordinates $(r, \theta)$ centered at $x$ and by (\ref{Lapdis}),
the Laplacian of the distance function $r$ satisfies
$$
\Delta r=\frac{\partial}{\partial r} \log \sigma(x, r, \theta).
$$
By the assumption, Ric$_{\infty} \geq -K$, $K\geq0$. According to Lemma \ref{Lapcompa}(2), we have
\begin{equation*}
\Delta r \leq\left.\frac{d}{d t}\left[\ln \chi_{0}(t)\right]\right|_{t=r},
\end{equation*}
where $\chi_{0}(t)=\left[s_{-K/(n-1)}(t)\right]^{n-1} e^{\alpha t}$. Thus we can get
\begin{equation*}
\Delta r \leq \alpha+(n-1) ct_{-K/(n-1)}(r).
\end{equation*}
Observe that $ct_{-K/(n-1)}(r) \leq \sqrt{\frac{K}{n-1}}\left(\sqrt{\frac{n-1}{K}}\frac{1}{r}+1\right)$ by an elementary argument. Consequently, one obtains
\begin{equation*}
\frac{\partial}{\partial r} \log \sigma(x, r, \theta) \leq \alpha+\sqrt{(n-1)K}\left(\sqrt{\frac{n-1}{K}}\frac{1}{r}+1\right),
\end{equation*}
which implies that
\begin{equation*}
\frac{\sigma\left(x, r_{2}, \theta\right)}{\sigma\left(x, r_{1}, \theta\right)} \leq\left(\frac{r_{2}}{r_{1}}\right)^{n-1} e^{r_{2}\left(\alpha+\sqrt{(n-1)K}\right)}.
\end{equation*}
This is just (\ref{ssigma}). Further, for any $ 0<s<r_{1}<t<r_{2}<R$, we have
\begin{equation*}
\sigma(x, t, \theta) s^{n-1} \leq t^{n-1} \sigma(x, s, \theta) e^{t\left(\alpha+\sqrt{(n-1)K}\right)} \leq t^{n-1} \sigma(x, s, \theta) e^{r_2\left(\alpha+\sqrt{(n-1)K}\right)} .
\end{equation*}
Now, integrating in $t$ from $r_{1}$ to $r_{2}$, we get
\begin{equation*}
s^{n-1} \int_{r_{1}}^{r_{2}} \sigma(x, t, \theta) d t \leq \frac{1}{n}\left(r_{2}^{n}-r_{1}^{n}\right) \sigma(x, s, \theta) e^{r_2\left(\alpha+\sqrt{(n-1)K}\right)} .
\end{equation*}
Then, integrating on both sides of above inequality with respect to $s$ from 0 to $r_{1}$ yields
\begin{equation*}
\frac{1}{n} r_{1}^{n} \int_{r_{1}}^{r_{2}} \sigma(x, t, \theta) d t \leq \frac{1}{n}\left(r_{2}^{n}-r_{1}^{n}\right) \int_{0}^{r_{1}} \sigma(x, s, \theta) e^{r_2\left(\alpha+\sqrt{(n-1)K}\right)} d s .
\end{equation*}
Further, we have
\begin{equation*}
\int_{r_{1}}^{r_{2}} d t \int_{\mathcal{D}_{x}(t)} \sigma(x, t, \theta) d \theta \leq \frac{r_{2}^{n}-r_{1}^{n}}{r_{1}^{n}} e^{r_2\left(\alpha+\sqrt{(n-1)K}\right)} \int_{0}^{r_{1}} d s \int_{\mathcal{D}_{x}(s)} \sigma(x, s, \theta) d \theta.
\end{equation*}
Then that is
\begin{equation*}
m\left(B_{r_{2}}\left(x\right)\right)-m\left(B_{r_{1}}\left(x\right)\right) \leq \frac{r_{2}^{n}-r_{1}^{n}}{r_{1}^{n}}  e^{r_2\left(\alpha+\sqrt{(n-1)K}\right)} m\left(B_{r_{1}}\left(x\right)\right) .
\end{equation*}
Therefore, it follows that
\begin{equation*}
\frac{m\left(B_{r_{2}}\left(x\right)\right)}{m\left(B_{r_{1}}\left(x\right)\right)}\leq\left(\frac{r_{2}}{r_{1}}\right)^{n} e^{r_2\left(\alpha+\sqrt{(n-1)K}\right)},
\end{equation*}
which is just (\ref{volBallS}). This completes the proof of Proposition \ref{volcom1}(1).

By Lemma \ref{Lapcompa}(3) and a similar argument as in (1), we can get Proposition \ref{volcom1}(2). \qed

\vskip 2mm

Similarly, we can derive the following volume comparison theorem by Lemma \ref{Lapcompa}.

\begin{prop}\label{volcom2}
Let $(M, F, m)$ be an n-dimensional Finsler manifold with a smooth volume form $dm$. Assume that ${\rm Ric}_{\infty} \geq K>0$. Then, along any minimizing geodesic starting from the center $x$ of $B_{R}(x)$, the following conclusions hold.
\ben
\item[(1)] If $\mathbf{S} \geq -\alpha$, $\alpha>0$, then, for any $0<r_1<r_2<\min \left\{R, \frac{\pi}{2} \sqrt{\frac{n-1}{K}}\right\}$,  we have
\be
\frac{\sigma\left(x, r_2, \theta\right)}{\sigma\left(x, r_1, \theta\right)} \leq\left(\frac{r_2}{r_1}\right)^{n-1} e^{r_2 \alpha}.
\ee
Further, we have
\be
\frac{m\left(B_{r_2}\left(x\right)\right)}{m\left(B_{r_1}\left(x\right)\right)} \leq\left(\frac{r_2}{r_1}\right)^n e^{r_2 \alpha}.
\ee

\item[(2)] If $|\tau| \leq k$, then, for any $0<r_1<r_2<\min \left\{R, \frac{\pi}{4} \sqrt{\frac{n-1}{K}}\right\}$,  we have
\be
\frac{\sigma\left(x, r_2, \theta\right)}{\sigma\left(x, r_1, \theta\right)} \leq\left(\frac{r_2}{r_1}\right)^{n+4 k-1}.
\ee
Further,  we have
\be
\frac{m\left(B_{r_2}\left(x\right)\right)}{m\left(B_{r_1}\left(x\right)\right)} \leq\left(\frac{r_2}{r_1}\right)^{n+4 k}.
\ee
\een
\end{prop}

Proposition \ref{volcom1} and Proposition \ref{volcom2} will play a very important role in following discussions.

\section{Local uniform Poincar\'{e} inequality and Sobolev inequality}\label{PSIne}

In view of the standard theory for general metric measure spaces, one needs the volume doubling condition and a local uniform Poincar\'{e} inequality to prove the local Sobolev inequality. As we known, C Xia proved the local uniform Poincar\'{e} inequality under the condition that ${\rm Ric}_N\geq -K$ for $N\in [n, \infty)$ and $K\geq 0$ \cite{CXia}. Based on the volume comparisons in Section \ref{volcomth}, we get the following local uniform Poincar\'{e} inequality by a similar argument.

\begin{thm}{\label{poineqd}}
Let $(M, F, m)$ be a forward complete Finsler manifold with finite reversibility $\Lambda$. Assume that ${\rm Ric}_{\infty} \geq K$ and $|\tau|\leq k$ for some $K \in \mathbb{R} $ and $k>0$. Then, there exist positive constants $c_{i}=c_{i}(n, k, \Lambda)(i=1,2)$ depending only on $n, k$ and the reversibility $\Lambda$ of $F$,  such that
\begin{equation}
\int_{B_{R}}\left|u-\bar{u}\right|^2 dm \leq c_{1} e^{c_{2} \sqrt{|K|}R} R^2 \int_{B_{R}} F^{*2}(du) dm  \label{pi12}
\end{equation}
for $u \in W_{\mathrm{loc}}^{1,2}(M)$, where $B_{R}=B_{R}^{+}(x_{0})$ is the forward geodesic ball of radius $R$ $(\leq \frac{\pi}{4}\sqrt{\frac{n-1}{K}}$ when $K> 0)$ for any $x_{0} \in M$ and $\bar{u}:=\frac{1}{m\left(B_{R}\right)} \int_{B_{R}}u~dm$.

In particular, if $ {\rm Ric}_{\infty} \geq 0$,  the above inequality is reduced to
\begin{equation}
\int_{B_{R}}\left|u-\bar{u}\right|^2 dm \leq c_1R^2 \int_{B_{R}} F^{*2}(du) dm. \label{p21}
\end{equation}
\end{thm}
{\it Proof.} For any $x, y \in B_{R}$, let $\eta:[0, d(x,y)] \rightarrow M$ be a unit speed minimizing geodesic of $F$ from $x$ to $y$, and $\bar{\eta}: [0, d(y,x)] \rightarrow M$ the unit speed minimizing geodesic of  $F$ from $y$ to $x$. Further, by using such geodesic $\eta$, we can define a map $\varphi_{x,t}: M \rightarrow M$ by $\varphi_{x,t}(y)= \eta (t)$ for any $y\in M$.

Observe that
\begin{equation*}
-\int_{0}^{d(y,x)} F(\nabla u(\bar{\eta}(t))) dt \leq u(y)-u(x) \leq \int_{0}^{d(x,y)} F(\nabla u(\eta(t))) dt
\end{equation*}
by (\ref{co-Finsler}). Then we have
\beq
\int_{B_{R}}\left|u-\bar{u}\right|^2 dm &\leq & \frac{1}{m\left(B_{R}\right)} \int_{B_{R}} \int_{B_{R}}|u(x)-u(y)|^2 dm(x) dm(y) \nonumber \\
& \leq & \frac{1}{m\left(B_{R}\right)} \int_{B_{R}} \int_{B_{R}}\left(\int_{0}^{d(x,y)} F(\nabla u(\eta(t))) dt \right.  \nonumber \\
&& \left.+\int_{0}^{d(y,x)} F(\nabla u(\bar{\eta}(t))) \mathrm{d} t\right)^2 dm(x) dm(y)  \nonumber\\
& \leq & \frac{4(\Lambda+1)R}{m\left(B_{R}\right)} \int_{B_{R}} \int_{B_{R}} \int_{0}^{d(x,y)} F^2(\nabla u(\eta(t))) dt dm(x) dm(y). \label{p1-1}
\eeq
Here, we have use $d(x,y)\leq d(x, x_{0})+d(x_{0},y)\leq (\Lambda+1)R$ in the last inequality. In the following, for simplicity of notation, we denote by $d=d(x,y)$.

Note that $\overleftarrow{\eta}(t):=\eta(d-t)$ is the unit speed minimizing geodesic from $y$ to $x$ with respect to the reversed metric $\overleftarrow{F}$. It is easy to see that $\overleftarrow{F}(\overleftarrow{\nabla} (-u))=F(\nabla(u))$ since $\overleftarrow{\nabla} (-u)=-\nabla u$. Thus
\be
\int_{0}^{\frac{d}{2}} F^2(\nabla u(\eta(t))) dt=\int_{\frac{d}{2}}^{d} F^2(\nabla u(\eta(d-t)))dt
=\int_{\frac{d}{2}}^{d} \overleftarrow{F}^2(\overleftarrow{\nabla}(-u)(\overleftarrow{\eta}(t)))dt. \label{p1-2}
\ee
Together (\ref{p1-2}) with (\ref{p1-1}), one obtains
\beq
\int_{B_{R}}\left|u-\bar{u}\right|^2 dm &\leq & \frac{4(\Lambda+1)R}{m\left(B_{R}\right)}  \int_{B_{R}} \int_{B_{R}} \int_{\frac{d}{2}}^{d} \big\{F^2(\nabla u(\eta(t)))  \nonumber\\
&& +\overleftarrow{F}^2(\overleftarrow{\nabla}(-u)(\overleftarrow{\eta}(t)))\big\} dt d m(x) dm(y). \label{p1-3}
\eeq
For any $z=\eta\left(t\right)$, we can decompose the measure along $\eta(t)$ by polar coordinates $(t, \theta)$ as
\be
dm(z)=e^{-\Psi(\eta(t))} {\rm Vol}_{g_{\dot{\eta}}}= \sigma(x, t, \theta) dt d\theta. \label{p1-4}
\ee
Therefore, for $t\in [d/2,d]$, we have by Proposition \ref{volcom1} and \ref{volcom2}
\begin{equation}{\label{psigma}}
\frac{\sigma(x, d, \theta)}{\sigma(x, t, \theta)} \leq 2^{n+4k-1} e^{(\Lambda+1)(n+4k-1)\frac{1}{\sqrt{n-1}}\sqrt{|K|}R}.
\end{equation}
Denote $c_n:=2^{n+4k-1}$ and $c_{n,\Lambda}:=(\Lambda+1)(n+4k-1)\frac{1}{\sqrt{n-1}}$. Then, by $(\ref{psigma})$, we get the following
\beq
&&\int_{B_{R}} \int_{B_{R}} \int_{\frac{d}{2}}^{d} F^2(\nabla u(\eta(t))) dt dm(x) dm(y) \nonumber\\
&&\leq c_n e^{c_{n,\Lambda}\sqrt{|K|}R} \int_{B_R} \int_{B_R} \int_{\frac{d}{2}}^{d} F^2(\nabla u(\eta(t))) \frac{\sigma(x, t, \theta)}{\sigma(x, d, \theta)} d t d m(y) d m(x)  \nonumber\\
&&\leq c_n e^{c_{n,\Lambda}\sqrt{|K|}R}\int_0^{(\Lambda+1) R} \int_{B_R} \int_{B_R} F^2\left(\nabla u(\eta(t))\right) \frac{\sigma(x, t, \theta)}{\sigma(x, d, \theta)} d m(y) d m(x) d t  \nonumber\\
&&\leq c_n e^{c_{n,\Lambda}\sqrt{|K|}R}\int_0^{(\Lambda+1) R} \int_{B_R} \int_{\varphi_{x,t}(B_R)} F^2\left(\nabla u(z)\right) dm(z) dm(x) dt  \nonumber\\
&&\leq c_n(\Lambda+1) e^{c_{n,\Lambda}\sqrt{|K|}R} m\left(B_R\right)R \int_{B_{(\Lambda+2) R}} F^2(z, \nabla u(z)) d m(z), \label{p1-5}
\eeq
where we have used $d\left(x, z\right)=t \leq d = d\left(x, y\right) \leq (\Lambda+1) R$ in the last inequality. Since the Ricci lower bound is the same for $F$ and $\overleftarrow{F}$, we get in the same way as above
\begin{equation}{\label{p1-6}}
\begin{aligned}
& \int_{B_R} \int_{B_R} \int_{\frac{d}{2}}^{d} \overleftarrow{F}^2(\overleftarrow{\nabla}(-u)(\overleftarrow{\eta}(t))) dt d m(x)d m(y) \\
& \quad\leq c_n(\Lambda+1) e^{c_{n,\Lambda}\sqrt{|K|}R} m\left(B_R\right) R\int_{B_{(\Lambda+2) R}} F^2(\nabla u(z))d m(z).
\end{aligned}
\end{equation}
Substituting ({\ref{p1-5}}) and ({\ref{p1-6}}) into ({\ref{p1-3}}), we conclude that
\begin{equation}{\label{p1-7}}
\int_{B_R}|u-\bar{u}|^2 d m \leq  \tilde{c}_1e^{\tilde{c}_2\sqrt{|K|}R} R^2  \int_{B_{(\Lambda+2) R}} F^2(z, \nabla u(z)) d m(z),
\end{equation}
where $\tilde{c}_1:=8c_n (\Lambda+1)^2$ and $\tilde{c}_2:=c_{n,\Lambda}$.

Finally, by the Whitny-type covering argument (Corollary 5.3.5 in \cite{PS}), (\ref{p1-7}) can be improved in a standard way to the inequality (\ref{pi12}). This completes the proof. \qed

\vskip 2mm
By a similar argument, we obtain the following
\begin{thm}{\label{poineqs}}
Let $(M, F, m)$ be a forward complete Finsler manifold with finite reversibility $\Lambda$. Assume that ${\rm Ric}_{\infty} \geq K$ and $S\geq -\alpha$ for some $K \in R$ and $\alpha>0$. Then, there exist positive constants $c_{i}=c_{i}(n, \Lambda)(i=1, 2)$ depending only on $n$ and the reversibility $\Lambda$ of $F$, such that
\begin{equation}
\int_{B_{R}}\left|u-\bar{u}\right|^2 dm \leq c_{1} e^{c_{2} \left(\alpha+\sqrt{|K|}\right) R} R^2 \int_{B_{R}} F^{*2}(du)dm   \label{pi22}
\end{equation}
for $u \in W_{\mathrm{loc}}^{1,2}(M)$, where $B_{R}=B_{R}^{+}\left(x_{0}\right)$ is the forward geodesic ball of radius $R$ $(\leq  \frac{\pi}{2} \sqrt{\frac{n-1}{K}}$ if $ K>0)$ for any $x_{0} \in M$.
\end{thm}

In \cite{MW1}, O. Munteanu and J. Wang proved the local uniform Sobolev inequality in smooth metric measure spaces $(M, g, e^{-f}dv)$ (see Lemma 3.2, \cite{MW1}). Their argument only used the structure of metric spaces and volume comparison and local uniform Poincar\'{e} inequality. Following their argument,  C. Xia gave directly the local uniform Sobolev inequality on forward complete Finsler manifolds by setting $A(R)=\sqrt{K}R$ in \cite{CXia} under the condition that ${\rm Ric}_N\geq -K$ for some $K>0$ and  $N\in [n, \infty)$, where $A(R):=\sup _{x \in B_{3R}(p)}|f|(x)$. Next we prove the following local uniform Sobolev inequalities based on the Poincar\'{e} inequality (\ref{pi12})  and volume comparisons theorems in Section \ref{volcomth}.

\begin{thm}{\label{soineqd}}
Let $(M, F, m)$ be an $n$-dimensional forward  complete Finsler manifold with finite reversibility $\Lambda$. Assume that ${\rm Ric}_{\infty} \geq K$ and $|\tau|\leq k$ for some $ K \in \mathbb{R}$ and $k>0$. Then, there exist positive constants $\nu(n, k)>2$ and  $c=c(n, k, \Lambda)$ depending only on $n$, $k$ and the reversibility $\Lambda$ of $F$, such that
\begin{equation}
\left(\int_{B_{R}}\left|u-\bar{u}\right|^{\frac{2\nu}{\nu-2}} dm\right)^{\frac{\nu-2}{\nu}} \leq e^{c(1+\sqrt{|K|}R)} m(B_R)^{-\frac{2}{\nu}}R^2 \int_{B_{R}} F^{*2}(du) dm  \label{p2}
\end{equation}
for $u \in W_{\mathrm{loc}}^{1,2}(M)$, where $B_{R}=B_{R}^{+}\left(x_{0}\right)$ is the forward geodesic ball of radius $R$ $(\leq  \frac{\pi}{4} \sqrt{\frac{n-1}{K}}$ if $K >0)$ for any $x_{0} \in M$ and $\bar{u}:=\frac{1}{m\left(B_{R}\right)} \int_{B_{R}} u~d m$, $\nu=4 (n+4k)-2$.
Further,
\begin{equation}
\left(\int_{B_{R}}\left|u\right|^{\frac{2\nu}{\nu-2}} dm\right)^{\frac{\nu-2}{\nu}} \leq e^{c(1+\sqrt{|K|}R)} m(B_R)^{-\frac{2}{\nu}}R^2 \int_{B_{R}} \left(F^{*2}(du) + R^{-2} u^2\right) dm. \label{p2-1}
\end{equation}
\end{thm}
{\it Proof.}  Because Finsler metrics are not usually reversible,  here we give a proof for this theorem by a similar argument as in \cite{MW1} for readers' convenience.

For any $y \in B_R$,  let $\gamma(t):[0,l]\rightarrow M$ be a unit speed minimizing geodesic from $x_0=\gamma(0)$ to $y=\gamma(l)$. Define a set of values $t_0:=0$, $t_i:=\sum_{j=1}^i \frac{R}{2^j}$ for $1\leq i< i_0$, $t_{i_0}:=l$. %where $i_0$ is the largest integer so that $t_{i_0}< l$.
Denote $y_0:=\gamma(t_0)=x_0$, $y_i:=\gamma\left(t_i\right)$ and $y_{i_0}:=\gamma\left(t_{i_0}\right)=y$. Besides, define $B_i:=B_{R_i}^+(y_i)$ for $i<i_0$ and $B_i:=B_{R_i}^+(y)$ for $i \geq i_0$, where $R_i:=\frac{ R}{2^{i+1}}$.
%Obviously, the set of forward geodesic balls $B_i$ covers $\gamma(t)$.

Let $u_{B_i}:=m\left(B_i\right)^{-1} \int_{B_i} u\, dm$. Then $\lim \limits_{i \rightarrow \infty} u_{B_i}=u(y)$. Denote $D_i:=B_{\tilde{R}_i}^+(z_i) $, where $\tilde{R}_i=\frac{R}{\Lambda2^{i+3}}$ and $z_i:=\gamma\left(\sum_{j=1}^{i+1} \frac{R}{2^j}-\frac{ R}{\Lambda2^{i+3}}\right)$ for $i<i_0$ and $z_i:=\gamma(l)=y$ for $i\geq i_0$. By the triangle inequality, it is easy to check that $D_i\subset B_i \cap B_{i+1}$.
Thus, we have
\beqn
&& \left|u_{B_0}-u(y)\right| \leq \sum_{i \geq 0}\left|u_{B_i}-u_{B_{i+1}}\right| \leq \sum_{i \geq 0}\left(\left|u_{B_i}-u_{D_i}\right|+\left|u_{D_i}-u_{B_{i+1}}\right|\right)\\
&&=\sum_{0 \leq i< i_0} \left(\left|u_{B_i}-u_{D_i}\right|+\left|u_{D_i}-u_{B_{i+1}}\right|\right)+\sum_{i\geq i_0} \left(\left|u_{B_i}-u_{D_i}\right|+\left|u_{D_i}-u_{B_{i+1}}\right|\right).
\eeqn
Furthermore,  for $0 \leq i< i_0$,  it is easy to see $B_i\subset B_{2\Lambda R/2^i}^+(z_i)$ and $B_{i+1}\subset B_{R/2^{i+1}}^+(z_i)$, and then,
$$
\frac{m(B_i)}{m(D_i)}\leq \frac{m(B_{2\Lambda R/2^i}^+(z_i))}{m(B^+_{\tilde{R}_{i}}(z_i))}\leq \tilde{c}_1e^{\tilde{c}_2 \sqrt{|K|}R}
$$
and
$$
\frac{m(B_{i+1})}{m(D_i)}\leq \frac{m(B_{ R/2^{i+1}}^{+}(z_i))}{m(B^{+}_{\tilde{R}_{i}}(z_i))}\leq \tilde{c}_1e^{\tilde{c}_2\sqrt{|K|}R},
$$
where $\tilde{c}_1:=(2^4\Lambda^2)^{n+4k}$ and $\tilde{c}_2:=2\Lambda\frac{n+4k}{\sqrt{n-1}}$, and we have used volume comparisons in Propositions \ref{volcom1} and \ref{volcom2}. Then
$$
\left|u_{B_i}-u_{D_i}\right| \leq m\left(D_i\right)^{-1} \int_{B_i}\left|u-u_{B_i}\right|dm \leq \tilde{c}_1 \mathrm{e}^{\tilde{c}_2 \sqrt{|K|}R} m\left(B_i\right)^{-1} \int_{B_i}\left|u-u_{B_i}\right|dm.
$$
Similarly, we have
$$
\left|u_{D_i}-u_{B_{i+1}}\right|\leq \tilde{c}_1 \mathrm{e}^{\tilde{c}_2 \sqrt{|K|}R} m\left(B_{i+1}\right)^{-1} \int_{B_{i+1}}\left|u-u_{B_{i+1}}\right|dm .
$$
For $i\geq i_0$, it is obvious that $D_i=B_{\tilde{R}_i}^+(y)\subset B_{R_{i+1}}^+(y)\subset B_{R_{i}}^+(y)$.
Hence we obtain
\beq
\left|u_{B_0}-u(y)\right| & \leq & 2\tilde{c}_1 \mathrm{e}^{\tilde{c}_2 \sqrt{|K|}R} \sum_{i \geq 0} m\left(B_i\right)^{-1} \int_{B_i}\left|u-u_{B_i}\right| dm \nonumber\\
& \leq & 2\tilde{c}_1 \mathrm{e}^{\tilde{c}_2 \sqrt{|K|}R} \sum_{i \geq 0}\left(m\left(B_i\right)^{-1} \int_{B_i}\left|u-u_{B_i}\right|^2dm\right)^{\frac{1}{2}}  \nonumber\\
& \leq & c_3 \mathrm{e}^{c_4 \sqrt{|K|}R} \sum_{i \geq 0} \frac{R}{2^{i+1}}\left(m\left(B_i\right)^{-1} \int_{B_i}F^{*2}(d u) dm\right)^{\frac{1}{2}},\label{s-1}
\eeq
where, in the last inequality we have used Poincar\'{e} inequality (\ref{pi12}) and $c_3:= 2 \tilde{c}_1\sqrt{c_1}$ and $c_4:= \tilde{c}_2+\frac{1}{2}c_2$.

Notice that
$$
\left|u_{B_0}-u(y)\right|=c R^{-\frac{1}{2}} \sum_{i \geq 0}\left(\frac{R}{2^{i+1}}\right)^{\frac{1}{2}}\left|u_{B_0}-u(y)\right|,
$$
where $c=\sqrt{2}-1$. Thus, for $R_i:=\frac{R}{2^{i+1}}$, combining the above equation with (\ref{s-1}) yields
$$
\sum_{i \geq 0}\left(R_i\right)^{\frac{1}{2}}\left|u_{B_0}-u(y)\right| \leq \frac{c_3}{c} {e}^{c_4 \sqrt{|K|}R} R^{\frac{1}{2}} \sum_{i \geq 0} R_i\left(m \left(B_i\right)^{-1} \int_{B_i}F^{*2}(d u)dm\right)^{\frac{1}{2}} .
$$
Hence, there exists an $i$ (depending on $y$ ) so that
$$
\left|u_{B_0}-u(y)\right|^2 \leq c_5 {e}^{c_6 \sqrt{|K|}R}\left(R R_i\right) m\left(B_i\right)^{-1} \int_{B_i} F^{*2}(du) dm,
$$
where $c_5:=\frac{c_3^2}{c^2}$ and $c_6:=2c_4$. Since $B_i \subset B_{3 \Lambda R_i}^+\left(y\right)$, it follows that, for each $y \in B_R$ there exists $r_y>0$ so that
\begin{equation}
\left|u_{B_0}-u(y)\right|^2 \leq c_5 \mathrm{e}^{c_6 \sqrt{|K|}R}\left(R r_y\right) m\left(B_{r_y}^+\left(y\right)\right)^{-1} \int_{B_{r_y}^+\left(y\right) \cap B_R}F^{*2}(du) dm .  \label{s-2}
\end{equation}

According to  Proposition \ref{volcom1} and \ref{volcom2}, we have
\begin{equation}\label{s-3}
\frac{m\left(B_R^+(x_0)\right)}{m\left(B_{r_y}^+\left(y\right)\right)} \leq \frac{m \left(B_{2\Lambda R}^+(y)\right)}{m\left(B_{r_y}^+\left(y\right)\right)} \leq \hat{c}_1 \mathrm{e}^{\hat{c}_2 \sqrt{|K|}R}\left(\frac{R}{r_y}\right)^{n+4k},
\end{equation}
where $\hat{c}_1:=(2\Lambda)^{n+4k}$ and $\hat{c}_2:=2\Lambda\frac{n+4k}{\sqrt{n-1}}$.
Solving $r_y$ from (\ref{s-3}) and plugging into (\ref{s-2}) yield
\beq
\left|u_{B_0}-u(y)\right|^2 & \leq & c_7 {e}^{c_8 \sqrt{|K|}R} R^{2} m\left(B_R\right)^{-\frac{1}{n+4k}} m\left(B_{r_y}\left(y\right)\right)^{\frac{1}{n+4k}-1} \nonumber\\
&& \times \int_{B_{r_y}^+\left(y\right) \cap B_R}F^{*2}(du) dm,  \label{s-4}
\eeq
where $c_7:=c_5\hat{c}_1^{\frac{1}{n+4k}}$ and $c_8:=c_6+\frac{\hat{c}_2}{n+4k}$.

We now define $A_t:=\left\{y \in B_R \mid \left|u_{B_0}-u(y)\right| \geq t\right\}$. Applying the Vitali covering lemma (Lemma 3.4 in \cite{Simon}), we can find a countable disjoint collection $\left\{B_{r_i}^{+}\right\}_{i \in I}$ of balls from $\left\{B_{r_y}^{+}\left(y\right) \mid y \in A_t\right\}$ such that for any $y \in A_t$, there exists $i \in I$ such that $B_{r_i}^+ \cap B_{r_y}^+(y) \neq \emptyset$ and $B_{r_y}^+(y) \subset B_{5\Lambda r_i}^+$. Then, by Propositions \ref{volcom1} and \ref{volcom2} and (\ref{s-4}), we have
\beqn
m\left(A_t\right)^{1-\frac{1}{n+4k}} & \leq & \bar{c}_1 {e}^{\bar{c}_2 \sqrt{|K|}R} \sum_{i \in I} m\left(B_{r_i}^+\right)^{1-\frac{1}{n+4k}} \\
& \leq & c_9 {e}^{c_{10} \sqrt{|K|}R} \frac{R^2}{t^2} m\left(B_R\right)^{-\frac{1}{n+4k}} \sum_{i \in I} \int_{B_{r_i}^+\cap B_R}F^{*2}(du) dm \\
& = & c_9 {e}^{c_{10} \sqrt{|K|}R} \frac{R^2}{t^2} m\left(B_R\right)^{-\frac{1}{n+4k}} \int_{B_R}F^{*2}(du)dm,
\eeqn
where $\bar{c}_1:=(5\Lambda)^{n+4k}$, $\bar{c}_2:=5\Lambda\frac{n+4k}{\sqrt{n-1}}$, $c_{9}:=\bar{c}_1c_7$ and $c_{10}:=\bar{c}_2+c_8$, and we have used the fact that $B_{r_y}^+(y) \subset B_{5\Lambda r_i}^+$ and volume comparisons in Propositions \ref{volcom1} and \ref{volcom2} in the first inequality. It may be rewritten as
\be
m\left(A_t\right) \leq t^{-\frac{2 (n+4k)}{n+4k-1}} \mathcal{B}, \label{s-5}
\ee
where
$$
\mathcal{B}:=c_{11} \mathrm{e}^{c_{12}\sqrt{|K|}R} R^{\frac{2 (n+4k)}{n+4k-1}} m\left(B_R\right)^{-\frac{1}{n+4k-1}}\left(\int_{B_R} F^{*2}(du)dm \right)^{\frac{n+4k}{n+4k-1}}
$$
and $c_{11}:=c_9^{\frac{n+4k}{n+4k-1}}$, $c_{12}:=\frac{n+4k}{n+4k-1}c_{10}$. Now, for any $\frac{2 (n+4k)}{n+4k-1}>q>2$, applying Cavalieri's principle (Theorem 14.10 in \cite{met}), we have
\beqn
\int_{B_R^+(x_0)}\left|u-u_{B_0}\right|^{q} dm & = & q \int_0^{\infty} t^{q-1} m\left(A_t\right) d t \\
& =& q \int_{0}^{T} t^{q-1} m\left(A_t\right) d t+q \int_T^{\infty} t^{q-1} m\left(A_t\right) d t \\
& \leq & T^{q} m\left(B_R\right)+\frac{q}{\frac{2 (n+4k)}{n+4k-1}-q} T^{q-\frac{2 (n+4k)}{n+4k-1}} \mathcal{B},
\eeqn
where we have used $m(A_t)\leq m(B_R)$ and (\ref{s-5}) in the last inequality.

Now, choosing $q=\frac{2 (n+4k)-1}{n+4k-1}$ and  $T:=\left(R^2 \frac{1}{m \left(B_R\right)} \int_{B_R}F^{*2}(du)dm\right)^{\frac{1}{2}}$, we get
$$
\left(\int_{B_R}\left|u-u_{B_0}\right|^{\frac{2 \nu}{\nu-2}}\right)^{\frac{\nu-2}{\nu}} \leq c_{13} \mathrm{e}^{c_{12} \sqrt{|K|}R} \frac{R^2}{m\left(B_R\right)^{\frac{2}{\nu}}} \int_{B_R}F^{*2}(du)dm,
$$
where $\nu:=\frac{2 q}{q-2}=4 (n+4k)-2$, $c_{13}:=2(n+4k)c_{11}$, from which Sobolev inequality (\ref{p2}) holds.

Moreover, by triangle inequality in $L^{p}\left(B_{R}\right)$ $(p= \frac{2\nu}{\nu -2})$ and $(a+b)^{2}\leq 2(a^{2}+b^{2})$ for $a, b \geq 0$, we have
\begin{equation*}
\left(\int_{B_{R}}\left|u\right|^{\frac{2\nu}{\nu-2}} dm\right)^{\frac{\nu-2}{\nu}}
\leq  2 \left(\int_{B_{R}}\left|u-\bar{u}\right|^{\frac{2\nu}{\nu-2}} dm\right)^{\frac{\nu-2}{\nu}}+2 \left(\int_{B_{R}}\left|\bar{u}\right|^{\frac{2\nu}{\nu-2}} dm\right)^{\frac{\nu-2}{\nu}},
\end{equation*}
from which and (\ref{p2}), we obtain (\ref{p2-1}). This completes the proof. \qed

\vskip 2mm

By a similar argument as in Theorem \ref{soineqd}, we have the following theorem.
\begin{thm}{\label{soineqs}}
Let $(M, F, m)$ be a forward complete Finsler manifold with finite reversibility $\Lambda$. Assume that ${\rm Ric}_{\infty} \geq K$ and $S\geq -\alpha$ for some $ K \in \mathbb{R}$ and  $\alpha>0$. Then, there exist positive constants $\nu> 2$ and $c=c(n, \Lambda)$ depending on $n$ and the reversibility $\Lambda$ of $F$,  such that
\begin{equation}
\left(\int_{B_{R}}\left|u-\bar{u}\right|^{\frac{2\nu}{\nu-2}} dm\right)^{\frac{\nu-2}{\nu}} \leq e^{c\left(1+ \left(\alpha+\sqrt{|K|}\right) R\right)} m(B_R)^{-\frac{2}{\nu}}R^2 \int_{B_{R}} F^{*2}(du) dm     \label{s1}
\end{equation}
for $u \in W_{\mathrm{loc}}^{1,2}(M)$, where $B_{R}=B_{R}^{+}\left(x_{0}\right)$ is the forward geodesic ball of radius $R (\leq  \frac{\pi}{2} \sqrt{\frac{n-1}{K}}$ if $K>0)$ for any $x_{0} \in M$ and $\bar{u}:=\frac{1}{m\left(B_{R}\right)} \int_{B_{R}} u~d m$. Further,
\begin{equation*}
\left(\int_{B_{R}}\left|u\right|^{\frac{2\nu}{\nu-2}} dm\right)^{\frac{\nu-2}{\nu}} \leq e^{c\left(1+ \left(\alpha+\sqrt{|K|}\right) R\right)} m(B_R)^{-\frac{2}{\nu}}R^2 \int_{B_{R}} \left(F^{*2}(du) + R^{-2} u^2\right) dm.
\end{equation*}
\end{thm}

\section{Proofs of the main theorems}\label{proof}

In order to prove the Harnack inequality and then the gradient estimate for harmonic functions, we need the mean value inequality for nonnegative subsolutions of a class of elliptic operators. Hence, we first give the proof of mean value inequality under the condition that ${\rm Ric}_\infty \geq K$ for $K \in \mathbb{R}$.
\vskip 2mm

{\it Proof of Theorem \ref{meanineq} } Since $u$ is a nonnegative function satisfying $\Delta u\geq-fu$ in the weak sense on $B_R$, we have
\begin{equation}{\label{mean-2}}
\int_{B_{R}} d \phi(\nabla u) d m \leq \int_{B_{R}}\phi fu dm
\end{equation}
for any nonnegative function $\phi \in \mathcal{C}_0^{\infty}\left(B_R\right)$.  For any $0 < \delta<\delta^{\prime} \leq 1$ and $a\geq1$, let $\phi=u^{2a-1}\varphi^2$,  where $\varphi$ is a cut-off function defined by
$$
\varphi(x)= \begin{cases}1 & \text { on } B_{\delta R}, \\ \frac{\delta^{\prime} R-d \left(x_0, x\right)}{\left(\delta^{\prime}-\delta\right) R} & \text { on } B_{\delta^{\prime} R} \backslash B_{\delta R}, \\ 0 & \text { on } M \backslash B_{\delta^{\prime} R}.\end{cases}
$$
Then $F^*(-d \varphi) \leq \frac{1}{\left(\delta^{\prime}-\delta\right) R}$ and hence $F^*(d \varphi) \leq \frac{\Lambda}{\left(\delta^{\prime}-\delta\right) R}$ a.e. on $B_{\delta^{\prime} R}$. Thus, by (\ref{mean-2}), we have
$$
(2a-1)\int_{B_R} \varphi^2 u^{2a-2}F^{*2}(d u) d m  +2 \int_{B_R} \varphi u^{2a-1} d \varphi(\nabla u ) d m \leq  \int_{B_R} \varphi^2u^{2a} f d m.
$$
Then
$$
a^2\int_{B_R} \varphi^2 u^{2a-2}F^{*2}(d u) d m \leq -2a \int_{B_R} \varphi u^{2a-1} d \varphi(\nabla u ) d m + a \int_{B_R} \varphi^2 u^{2a} f d m,
$$
namely,
\beqn
\int_{B_R} \varphi^2 F^{*2}(d u^{a}) d m &\leq & -2 \int_{B_R} \varphi u^{a} d \varphi(\nabla u^{a} ) d m + a \int_{B_R} \varphi^2 u^{2a} f d m \\
&\leq & 2 \int_{B_R} \varphi u^{a} F^{*}(-d \varphi) F\left(\nabla u^{a}\right)dm + a\int_{B_R}\varphi^2 u^{2a} f d m \\
& \leq & \frac{1}{2} \int_{B_R} \varphi^2 F^2(\nabla u^{a}) dm + 2\int_{B_R}u^{2a} F^{*2}(-d \varphi) dm + a \int_{B_R} \varphi^2 u^{2a} f d m.\\
\eeqn
It follows from the above inequality that
\beq
\int_{B_R} \varphi^2 F^{*2}(d u^a) dm &\leq &  4 \int_{B_R} u^{2a} F^{* 2}(-d \varphi) d m + 2a \int_{B_R} \varphi^2 u^{2a} f d m \nonumber \\
&\leq & \frac{4}{\left(\delta^{\prime}-\delta\right)^2 R^2}\int_{B_{\delta^{'}R}} u^{2a} dm + 2a \mathcal{A}\int_{B_{\delta^{'}R}} u^{2a} d m, \label{meanv1}
\eeq
where $\mathcal{A}:=\sup\limits_{ B_{R}}|f|$. Moreover,
$$
\int_{B_R} \varphi^{2a} F^{*2}(d u^a) dm \leq \int_{B_R} \varphi^{2} F^{*2}(d u^a) dm
$$
and
$$
\int_{B_R} u^{2a} F^{*2}(d \varphi^a) dm =a^2\int_{B_R} u^{2a} \varphi^{2a-2} F^{*2}(d \varphi) dm
\leq  \frac{a^2\Lambda^2}{(\delta'-\delta)^2R^2} \int_{B_{\delta'R}} u^{2a} dm.
$$
In sum, by H\"{o}lder's inequality and Sobolev inequality (\ref{p2-1}), we have
\beqn
&& \int_{B_{\delta R}} u^{2a\left(1+\frac{2}{\nu}\right)}d m \leq\int_{B_R}(u \varphi)^{2a\left(1+\frac{2}{\nu}\right)} d m \leq\left(\int_{B_R}(u^a \varphi^a)^{\frac{2 \nu}{\nu-2}} d m\right)^{\frac{\nu-2}{\nu}} \cdot\left(\int_{B_R}(u^a \varphi^a)^2 d m\right)^{\frac{2}{\nu}} \\
&& \leq \mathcal{B} \int_{B_R}\left(F^{* 2}(d(u^a \varphi^a))+R^{-2} u^{2a} \varphi^{2a}\right) d m \cdot\left(\int_{B_R} u^{2a} \varphi^{2a} d m\right)^{\frac{2}{\nu}} \\
&& \leq \mathcal{B} \int_{B_R}\left(2 \varphi^{2a} F^{* 2}(d u^a)+2 u^{2a} F^{* 2}(d \varphi^a)+R^{-2} u^{2a} \varphi^{2a}\right) d m \cdot\left(\int_{B_{\delta^{\prime} R}} u^{2a} d m\right)^{\frac{2}{\nu}} \\
&& \leq \mathcal{B}\left(\frac{8+2a^2\Lambda^2}{\left(\delta^{\prime}-\delta\right)^2 R^2}+\frac{1}{R^2}+4a\mathcal{A}\right)\left(\int_{B_{\delta^{\prime} R}} u^{2a} d m\right)^{1+\frac{2}{\nu}}\\
&& \leq 11a^2\Lambda^2\mathcal{B}\left(\frac{1+\mathcal{A}R^2}{\left(\delta^{\prime}-\delta\right)^2 R^2}\right)\left(\int_{B_{\delta^{\prime} R}} u^{2a} d m\right)^{1+\frac{2}{\nu}},
\eeqn
where $\mathcal{B}:=e^{c(1+\sqrt{|K|} R)} R^2 m\left(B_R\right)^{-2 / \nu}$, $\nu$ and $c$ were chosen as in Theorem \ref{soineqd}, and we have used the fact that $F^{*2}(du+d\varphi)\leq 2F^{*2}(du)+2F^{*2}(d\varphi)$ in the third line. Let $t:=1+\frac{2}{\nu}$, the above inequality implies that
\begin{equation}\label{mean-3}
\int_{B_{\delta R}} u^{2 a t } d m \leq \frac{\mathcal{B}\Theta}{\left(\delta^{\prime}-\delta\right)^2 R^2}\left(\int_{B_{\delta^{\prime} R}} u^{2 a } d m\right)^t,
\end{equation}
where $\Theta:=11a^2\Lambda^2(1+\mathcal{A}R^2)$. For any $0<\delta <1$, let $\delta_0=1$ and $\delta_{i+1}=\delta_i-\frac{1-\delta}{2^{i+1}}$ on $B_{\delta_i R}$, $i=0,1, \cdots$. Applying (\ref{mean-3}) for $\delta^{\prime}=\delta_i, \delta=\delta_{i+1}$ and $a=t^i$, we have
$$
\int_{B_{\delta_{i+1} R}} u^{2 t^{i+1}} d m \leq \frac{4^{i+1}\mathcal{B}\Theta}{((1-\delta) R)^2}\left(\int_{B_{\delta_i R}} u^{2 t^i} d m\right)^t.
$$
By iteration, one obtains that
\beqn
\|u^2\|_{L^{t^{i+1}}(B_{\delta_{i+1} R})}&=& \left(\int_{B_{\delta_{i+1} R}} u^{2 t^{i+1}} d m\right)^{\frac{1}{t^{i+1}}} \\
&\leq & 4^{\sum j t^{-j}}\left(\mathcal{B}\Theta)\right)^{\sum t^{-j}}[(1-\delta) R]^{-2 \sum t^{-j}} \cdot \int_{B_R} u^2 d m,
\eeqn
in which $\sum$ denotes the summation on $j$ from 1 to $i+1$. Since $\sum_{j=1}^{\infty} t^{-j}=\frac{\nu}{2}$ and $\sum_{j=1}^{\infty} j t^{-j}=\frac{\nu^2+2\nu}{4}$, we have
\beq
\|u^2\|_{L^{\infty} (B_{\delta R})} &\leq & c_0 \mathcal{B}^{\frac{\nu}{2}}\Theta^{\frac{\nu}{2}}(1-\delta)^{-\nu}R^{-\nu}\int_{B_R} u^2 d m \nonumber\\
&= & e^{\tilde{C}(1+\sqrt{|K|} R)}(1+\mathcal{A}R^2)^{\frac{\nu}{2}} (1-\delta)^{-\nu} m\left(B_R\right)^{-1}  \int_{B_R} u^2 d m,\label{meanp2}
\eeq
which implies (\ref{meanineq-1}) with $p=2$, where $c_0=4^{\frac{\nu^2+2\nu}{4}}$ and $\tilde{C}:=\log(c_0(11a^2\Lambda^2)^{\frac{\nu}{2}})+\frac{\nu}{2}c>0$. This finishes the proof in the case when $p=2$.

Next we will consider the case when $0<p<2$ by using Moser iteration again. For any $0<\delta<1$, choose $\varepsilon\in(0,1)$ with $0<\delta+\varepsilon\leq1$. Then (\ref{meanp2}) implies
\beq
\sup\limits_{B_{\delta R}} u^{2}
&\leq& e^{\tilde{C}(1+\sqrt{|K|}R)}(1+\mathcal{A}R^2)^{\frac{\nu}{2}}\varepsilon^{-\nu}m(B_{R})^{-1}\int_{B_{(\delta+\varepsilon)R}} u^2 dm \nonumber\\
&\leq& e^{\tilde{C}(1+\sqrt{|K|}R)}(1+\mathcal{A}R^2)^{\frac{\nu}{2}}\varepsilon^{-\nu} m(B_{R})^{-1}\left(\sup\limits_{B_{(\delta+\varepsilon)R}}u^2\right)^{1-\frac{p}{2}}\int_{B_{R}} u^p dm. \label{meanpp-1}
\eeq

Let $\lambda=1-\frac{p}{2}>0$ and $A(\delta):=\sup\limits_{B_{\delta R}}u^2$. Choose $\delta_0=\delta$ and $\delta_i=\delta_{i-1}+\frac{1-\delta}{2^i}$, $i=1,2\cdots$. Applying (\ref{meanpp-1}) for $\delta=\delta_{i-1}$ and $\delta+\varepsilon=\delta_{i}$, we have
$$
A(\delta_{i-1})\leq \mathcal{\tilde{B}}\ 2^{i\nu}(1-\delta)^{-\nu}A(\delta_{i})^{\lambda},
$$
where $\mathcal{\tilde{B}}:=e^{\tilde{C}(1+\sqrt{|K|}R)}(1+\mathcal{A}R^2)^{\frac{\nu}{2}}m(B_{R})^{-1}\int_{B_{R}} u^p dm.$
By iterating, we get
$$
A(\delta_0)\leq \mathcal{\tilde{B}}^{\sum_{i=1}^j\lambda^{i-1}} 2^{\nu\sum_{i=1}^ji\lambda^{i-1}}(1-\delta)^{-\nu\sum_{i=1}^j\lambda^{i-1}}A(\delta_{j})^{\lambda^j}.
$$
Since $\lim\limits_{j\rightarrow\infty}\delta_j=1$, $\lim\limits_{j\rightarrow\infty}\lambda^j=0$, $\sum_{i=1}^{\infty}\lambda^{i-1}=\frac{2}{p}$ and $\sum_{i=1}^{\infty}i\lambda^{i-1}$ converges, then we obtain by letting $j\rightarrow\infty$
$$
\sup\limits_{B_{\delta R}}u^p \leq A(\delta)^{\frac{p}{2}} \leq e^{C(1+\sqrt{|K|}R)}(1+\mathcal{A}R^2)^{\frac{\nu}{2}}(1-\delta)^{-\nu}m(B_{R})^{-1}\int_{B_{R}} u^p dm,
$$
where $C=C(n, k, \nu, p, \Lambda)>0$. This finishes the proof of Theorem \ref{meanineq}.
\qed

\vskip 2mm

From Theorem \ref{meanineq}, we can prove the following lemma, which is important for the subsequent discussions.

\vskip 2mm
\begin{lem}\label{Delta} \  Let $(M, F, m)$ be an $n$-dimensional forward complete Finsler measure space equipped with a uniformly convex and uniformly smooth Finsler metric $F$. Assume that ${\rm Ric}_{\infty}\geq -K$ and $|\tau|\leq k$ for some $K\geq 0$ and $k>0$. Suppose that $u \in W^{2,2}(B_R)$  satisfies
\begin{equation}\label{hF-1}
\int_{B_R}\phi\Delta^{\nabla u} F^2(\nabla u)dm\geq -2 K \int_{B_R} \phi F^2(\nabla u) dm
\end{equation}
for all nonnegative functions $\phi \in \mathcal{C}_0^{\infty}\left(B_R\right)$. Then for any $p \in(0,2]$ and $0<\delta <1$, there are constants $\nu>2$ and $C=C(n, k, \nu, p, \kappa, \kappa^*)$ such that
\begin{equation}\label{hF-2}
\sup _{B_{\delta R}} F(\nabla u)^{2p} \leq e^{C(1+\sqrt{K} R)}(1+2KR^2)^{\frac{\nu}{2}} (1-\delta)^{-\nu} m\left(B_R\right)^{-1}  \int_{B_R} F(\nabla u)^{2p} dm.
\end{equation}
\end{lem}
{\it Proof.} Let $h(x):=F^{2}(\nabla u)$. Then $h\in W^{1,2}(B_R)$. For any $0 < \delta<\delta^{\prime} \leq 1$ and $a\geq1$, let $\phi=h^{2a-1}\varphi^2$,  where $\varphi$ is a cut-off function defined by
$$
\varphi(x)= \begin{cases}1 & \text { on } B_{\delta R}, \\ \frac{\delta^{\prime} R-d \left(x_0, x\right)}{\left(\delta^{\prime}-\delta\right) R} & \text { on } B_{\delta^{\prime} R} \backslash B_{\delta R}, \\ 0 & \text { on } M \backslash B_{\delta^{\prime} R}.\end{cases}
$$
Then $F^*(-d \varphi) \leq \frac{1}{\left(\delta^{\prime}-\delta\right) R}$ and hence $F^*(d \varphi) \leq \frac{\kappa}{\left(\delta^{\prime}-\delta\right) R}$ a.e. on $B_{\delta^{\prime} R}$. Thus, by (\ref{hF-1}), we have
$$
(2a-1)\int_{B_R} \varphi^2 h^{2a-2}dh(\nabla^{\nabla u} h) d m  +2 \int_{B_R} \varphi h^{2a-1} d \varphi(\nabla^{\nabla u} h) d m \leq  2K\int_{B_R} \varphi^2h^{2a} d m.
$$
It follows from (\ref{unisc}) that
$$
a^2\tilde{\kappa}^*\int_{B_R} \varphi^2 h^{2a-2}F^{*2}(d h) d m \leq -2a \int_{B_R} \varphi h^{2a-1} d \varphi(\nabla^{\nabla u} h ) d m + 2aK \int_{B_R} \varphi^2 h^{2a} d m,
$$
namely,
\beqn
\tilde{\kappa}^*\int_{B_R} \varphi^2 F^{*2}(d h^{a}) d m &\leq & -2 \int_{B_R} \varphi h^{a} d \varphi(\nabla^{\nabla u} h^{a} ) d m + 2aK \int_{B_R} \varphi^2 h^{2a} d m \\
&\leq &2 \tilde{\kappa}\int_{B_R} \varphi h^{a} F^{*}(-d \varphi) F\left(d h^{a}\right)dm + 2aK\int_{B_R}\varphi^2 h^{2a} d m \\
& \leq & \frac{1}{2} \tilde{\kappa}^*\int_{B_R} \varphi^2 F^{*2}(d h^{a}) dm + \frac{2\tilde{\kappa}^2}{\tilde{\kappa}^*}\int_{B_R}h^{2a} F^{*2}(-d \varphi) dm\\ & &+ 2aK \int_{B_R} \varphi^2 h^{2a} d m,\\
\eeqn
where we have used $-d\varphi(\nabla^{\nabla u} h)\leq \tilde{\kappa}F^{*}(-d \varphi)F^{*}(d h)$. Then it follows from the above inequality that
\beqn
\int_{B_R} \varphi^2 F^{*2}(d h^a) dm &\leq &  \frac{4\tilde{\kappa}^2}{\tilde{\kappa}^{*2}} \int_{B_R} h^{2a} F^{* 2}(-d \varphi) d m + \frac{4a}{\tilde{\kappa}^{*}} K \int_{B_R} \varphi^2 h^{2a}  d m  \\
&\leq & \frac{4\tilde{\kappa}^2}{\tilde{\kappa}^{*2}\left(\delta^{\prime}-\delta\right)^2 R^2}\int_{B_{\delta^{'}R}} h^{2a} dm + \frac{4a}{\tilde{\kappa}^{*}} K \int_{B_{\delta^{'}R}} h^{2a} d m ,
\eeqn
which is an analogue of (\ref{meanv1}). And then, along the proof of Theorem \ref{meanineq} we can get (\ref{hF-2}). \qed

\vskip 2mm
Further, by Lemma \ref{Delta}, we can prove a Liouville property for positive harmonic functions.

\vskip 2mm
{\it Proof of Corollary \ref{Liouville}.} It follows from the Bochner-Weitzenb\"{o}ck type formula (\ref{BWforinf}) and Ric$_{\infty}\geq0$ that
\begin{equation*}
\int_M\phi\Delta^{\nabla u} F^2(\nabla u)dm\geq0
\end{equation*}
for all nonnegative functions $\phi \in \mathcal{C}_0^{\infty}\left(M\right)$. By Lemma \ref{Delta}, we have
\begin{equation}\label{g-11}
\sup\limits_{B_{\frac{1}{4} R}}F^2(\nabla u)\leq \frac{C}{m(B_{\frac{1}{2}R})}\int_{B_{\frac{1}{2}R}}F^2(\nabla u)dm.
\end{equation}
Let $\phi=u\varphi^2$, where $\varphi$ is a cut-off function defined by
$$
\varphi(x)= \begin{cases}1 & \text { on } B_{\frac{1}{2} R}, \\ \frac{ 2R-2d_F\left(x_0, x\right)}{ R} & \text { on } B_{R} \backslash B_{\frac{1}{2}R}, \\ 0 & \text { on } M \backslash B_{R}.\end{cases}
$$
Then $F^*(-d \varphi) \leq \frac{2}{ R}$ and hence $F^*(d \varphi) \leq \frac{2\kappa}{R}$ a.e. on $B_{R}$. Thus we have
\beqn
\int_{M}\varphi^2 F^2(\nabla u) dm &=&-2\int_{M}\varphi u d\varphi(\nabla u) dm\\
&\leq & \frac{1}{2}\int_{M}\varphi^2 F^2(\nabla u) dm +2\int_{M}u^2 F^{*2}(-d\varphi) dm.
\eeqn
Then it follows that
$$
\int_{B_{\frac{1}{2} R}}F^2(\nabla u) dm\leq \frac{16}{R^2}m\left({B_{R}}\right)\left(\sup\limits_{B_{R}} u\right)^2.
$$
Combining (\ref{g-11}) with (\ref{L-1}) yields
\begin{equation*}%\label{g-2}
\sup\limits_{B_{\frac{1}{4} R}}F^2(\nabla u)\leq \frac{C}{R^2}\left(\sup\limits_{B_{ R}}u\right)^2\rightarrow 0\;\; (R\rightarrow \infty),
\end{equation*}
from which, we conclude that $F(\nabla u)=0$ on $M$. Hence $u$ is a constant on $M$. \qed

\vskip 2mm

In order to prove Theorem \ref{mean2}, we need the following lemmas. First, we have the following lemma.

\begin{lem}{\rm(\cite{PS})}\label{measure} Suppose that $\{U_{\sigma} \mid 0< \sigma \leq 1\}$ is a family of measurable subsets of a measurable set $U \subset \mathbb{R}^{n}$ with the property that $U_{\sigma'}\subset U_{\sigma}$ if $\sigma'\leq\sigma$. Fix $0<\delta <1$. Let $\gamma$ and $C$ be positive constants and $0<\alpha_0\leq\infty$. Let $g$ be a positive measurable function defined on $U_1=U$ which satisfies
\begin{equation}\label{alpha}
\left(\int_{U_{\sigma'}}g^{\alpha_0}dm\right)^{\frac{1}{\alpha_{0}}}\leq \left[C(\sigma-\sigma')^{-\gamma}m(U)^{-1}\right]^{\frac{1}{\alpha}-\frac{1}{\alpha_0}} \left(\int_{U_{\sigma}}g^{\alpha}dm\right)^{\frac{1}{\alpha}}
\end{equation}
for all $\sigma$, $\sigma'$, $\alpha$ satisfying $0<\delta\leq\sigma'<\sigma\leq1$ and $0<\alpha\leq \min\{1, \frac{\alpha_0}{2}\}$. Assume further that $g$ satisfies
\be
m(\log g>\lambda)\leq Cm(U)\lambda^{-1} \label{Lcondi2}
\ee
for all $\lambda>0$. Then
\begin{equation}\label{U}
\left(\int_{U_{\delta}}g^{\alpha_0}dm\right)^{\frac{1}{\alpha_{0}}}\leq C_0m(U)^{\frac{1}{\alpha_0}},
\end{equation}
where $C_0$ depends only on $\delta$, $\gamma$, $C$ and a lower bound on $\alpha_0$.
\end{lem}

Further, we can prove the following lemma.

\begin{lem}\label{supu-1}
Let $(M, F, m)$ be an $n$-dimensional forward complete Finsler measure space with finite reversibility $\Lambda$. Assume that ${\rm Ric}_{\infty} \geq K$ and $|\tau|\leq k$ for some $K\in \mathbb{R}$  and $k>0$. If $u$ is a positive superharmonic function on $B_R$, then,  for any $\delta \in (0,1)$,  there exist positive constant $C=C\left(n, \delta, k, \Lambda\right)$ depending on $n, \delta, k$ and  $\Lambda$, such that
$$
\sup _{B_{\delta R}} u^{-1} \leq e^{C(1+\sqrt{|K|} R)} (1-\delta)^{-\nu} m\left(B_R\right)^{-1}  \int_{B_R} u^{-1} dm. \label{meansup}
$$
\end{lem}
{\it Proof. } Since $u$ is a positive superharmonic function, we have
\begin{equation}{\label{meanu-2}}
\int_{B_{R}} d \phi(\nabla u) d m\geq0
\end{equation}
for any function $\phi \in \mathcal{C}_0^{\infty}\left(B_R\right)$.  For any $0<\delta<\delta^{\prime} \leq 1$, let $\phi=-bu^{b-1}\varphi^2$ and $b<-1$, where $\varphi$ is a cut-off function defined by
$$
\varphi(x)= \begin{cases}1 & \text { on } B_{\delta R}, \\ \frac{\delta^{\prime} R-d \left(x_0, x\right)}{\left(\delta^{\prime}-\delta\right) R} & \text { on } B_{\delta^{\prime} R} \backslash B_{\delta R}, \\ 0 & \text { on } M \backslash B_{\delta^{\prime} R}.\end{cases}
$$
Then $F^*(-d \varphi) \leq \frac{1}{\left(\delta^{\prime}-\delta\right) R}$ and hence $F^*(d \varphi) \leq \frac{\Lambda}{\left(\delta^{\prime}-\delta\right) R}$ a.e. on $B_{\delta^{\prime} R}$. Thus, by (\ref{meanu-2}), we have
$$
b(b-1)\int_{B_R} \varphi^2 u^{b-2}F^{*2}(d u) d m  +2 b \int_{B_R} \varphi u^{b-1} d \varphi(\nabla u ) d m \leq 0.
$$
Set $w:=u^{\frac{b}{2}}$. Then $dw=-\left|\frac{b}{2}\right|u^{\frac{b}{2}-1} du$. Hence, we have
$$\frac{b^2}{4}u^{b-2}\Lambda^{-2} F^{*2}(d u)\leq F^{*2}(d w)=\frac{b^2}{4}u^{b-2}F^{*2}(-d u)\leq \frac{b^2}{4}u^{b-2}\Lambda^2 F^{*2}(d u).$$
Then we get
\beqn
4\int_{B_R} \varphi^2F^{*2}(d w) d m &\leq &\frac{4b(b-1)}{b^2}\int_{B_R} \varphi^2F^{*2}(d w) d m \leq -2 b \Lambda^2 \int_{B_R} \varphi u^{b-1} d \varphi(\nabla u ) d m\\
&\leq & \frac{1}{2\Lambda^2}b^2\int_{B_R} \varphi^2 u^{b-2} F^2(\nabla u) dm + 2 \Lambda^6 \int_{B_R} w^2 F^{*2}(d \varphi) dm\\
&\leq & 2 \int_{B_R} \varphi^2 F^2(\nabla w) dm + 2 \Lambda^6 \int_{B_R} w^2 F^{*2}(d \varphi) dm,
\eeqn
namely,
$$
\int_{B_R} \varphi^2F^{*2}(d w) d m \leq  \Lambda^6 \int_{B_R} w^2 F^{*2}(d \varphi) dm \leq \frac{\Lambda^8}{\left(\delta^{\prime}-\delta\right)^2 R^2} \int_{B_{\delta'R}} w^2 dm.
$$
Similar to the proof of Theorem \ref{meanineq}, by H\"{o}lder's inequality and Sobolev inequality (\ref{p2-1}), we have
\beqn
&& \int_{B_{\delta R}} w^{2\left(1+\frac{2}{\nu}\right)}d m \leq\int_{B_R}(w \varphi)^{2\left(1+\frac{2}{\nu}\right)} d m \leq\left(\int_{B_R}(w \varphi)^{\frac{2 \nu}{\nu-2}} d m\right)^{\frac{\nu-2}{\nu}} \cdot\left(\int_{B_R}(w \varphi)^2 d m\right)^{\frac{2}{\nu}} \\
&& \leq \mathcal{B} \int_{B_R}\left(2 \varphi^{2} F^{* 2}(d w)+2 w^{2} F^{* 2}(d \varphi)+R^{-2} w^{2} \varphi^{2}\right) d m \cdot\left(\int_{B_{\delta^{\prime} R}} w^{2} d m\right)^{\frac{2}{\nu}} \\
&& \leq \frac{5\Lambda^8\mathcal{B}}{\left(\delta^{\prime}-\delta\right)^2 R^2}\left(\int_{B_{\delta^{\prime} R}} w^{2} d m\right)^{1+\frac{2}{\nu}},
\eeqn
where $\mathcal{B}:=e^{c(1+\sqrt{|K|} R)} R^2 m\left(B_R\right)^{-2 / \nu}$, $\nu$ and $c$ were chosen as in Theorem \ref{soineqd}. Let $\tau:=1+\frac{2}{\nu}$, $b:=-\beta$. Then the above inequality implies that
$$
\int_{B_{\delta R}} (u^{-1})^{\beta\tau}d m \leq \frac{5\Lambda^8\mathcal{B}}{\left(\delta^{\prime}-\delta\right)^2 R^2}\left(\int_{B_{\delta^{\prime} R}} (u^{-1})^{\beta} d m\right)^{\tau}.
$$
By iteration, one obtains that
$$
\sup _{B_{\delta R}} u^{-1} \leq e^{C(1+\sqrt{|K|} R)} (1-\delta)^{-\nu} m\left(B_R\right)^{-1}  \int_{B_R} u^{-1} dm.
$$
This finishes the proof. \qed

\vskip 2mm
In the following, we will prove Harnack inequality for harmonic functions.
\vskip 2mm

{\it Proof of Theorem \ref{mean2} }\,
By $L^p$ $(0<p\leq2)$  mean value inequality in Theorem \ref{meanineq} with $f=0$, for $\rho\in(0,1)$ and $0<\delta<1$, we have
\be
\sup\limits_{B_{\rho \delta R}}u \leq e^{C(1+\sqrt{|K|} R)} (1-\rho)^{-\nu} m\left(B_{\delta R}\right)^{-1}  \int_{B_{\delta R}} u dm. \label{mean2-1}
\ee

Next, we will combine the above inequality and Theorem \ref{poineqd} to prove
\begin{equation}
\int_{B_{\delta R}} u dm \leq e^{C_0(1+ \sqrt{|K|}R)}m(B_{R})\inf\limits_{B_{\delta R}}u     \label{mean2-2}
\end{equation}
for any $0<\delta<\delta'\leq1$, where $C_0:=C(n, k, \nu, p, \Lambda)>0$.

Since $u$ is a positive harmonic function, we have
$$
\int_{B_R} d\phi(\nabla u) dm=0
$$
for all nonnegative functions $\phi \in \mathcal{C}_0^{\infty}\left(B_R\right)$.   For any $0<\delta<\delta'\leq1$,  let $\phi=u^{-1}\varphi^2$, where $\varphi$ is a cut-off function defined by
$$
\varphi(x)= \begin{cases}1 & \text { on } B_{\delta R}, \\ \frac{\delta'R-d \left(x_0, x\right)}{(\delta'-\delta)R} & \text { on } B_{\delta' R} \backslash B_{\delta R}, \\ 0 & \text { on } M \backslash B_{\delta'R}.\end{cases}
$$
Then $F^*(-d \varphi) \leq \frac{1}{(\delta'-\delta)R}$ and hence $F^*(d \varphi) \leq \frac{\Lambda}{(\delta'-\delta)R}$ a.e. on $B_{\delta'R}$. Thus,
\beqn
0 &=&  \int_{B_R}d(u^{-1}\varphi^2)(\nabla u) dm = -\int_{B_R}u^{-2}\varphi^2 F^2(\nabla u)dm+ 2\int_{B_R} u^{-1}\varphi d\varphi(\nabla u)dm \\
&=&-\int_{B_R}\varphi^2F^2(\nabla v)dm+2\int_{B_R}\varphi d\varphi(\nabla v)dm,
\eeqn
where $v=\log u$. Then
\beqn
\int_{B_R}\varphi^{2}F^{2}(\nabla v)dm &=& 2\int_{B_R}\varphi d\varphi(\nabla v)dm \leq 2\int_{B_R}\varphi F(\nabla v)F^{*}(d \varphi) dm  \\
&\leq & \frac{1}{2}\int_{B_R} \varphi^2F^2(\nabla v)dm+2\int_{B_R}F^{*2}(d\varphi)dm.
\eeqn
Therefore,
\be
\int_{B_{\delta R}}F^2(\nabla v)dm\leq4\int_{B_{\delta' R}}F^{*2} (d\varphi)dm\leq\frac{4\Lambda^2}{(\delta'-\delta)^2R^2}m({B_{\delta'R}}).  \label{h1}
\ee
In addition, it follows that
\beq
t \ m\left(B_{\delta R}\bigcap \{|v-\bar{v}|\geq t\}\right) & \leq &\int_{B_{\delta R}}|v-\bar{v}|dm \nonumber\\
&\leq &\left(\int_{B_{\delta R}}|v-\bar{v}|^2dm\right)^{\frac{1}{2}}m(B_{\delta R})^{\frac{1}{2}} \nonumber\\
&\leq &\,\sqrt{c_1}e^{\frac{1}{2}c_2\sqrt{|K|}R}R\left(\int_{B_{\delta R}}F^2(\nabla v)dm\right)^{\frac{1}{2}}m(B_{\delta R})^{\frac{1}{2}}\nonumber\\
& \leq &\,\sqrt{c_1}e^{\frac{1}{2}c_2\sqrt{|K|}R} \frac{2\Lambda}{\delta'-\delta}  m(B_{R}),  \label{h2}
\eeq
where $\bar{v}=\frac{1}{m\left(B_{\delta R}\right)} \int_{B_{\delta R}} v dm$ and we have used Poincar\'{e} inequality (\ref{pi12}) in the third inequality and the last inequality follows from (\ref{h1}).

Letting $g=e^{-\bar{v}}u$ in Lemma \ref{measure}, we have $\Delta^{\nabla u} g=e^{-\bar{v}} \Delta u=0$. By the mean value inequality in Theorem \ref{meanineq} with $p=\frac{1}{2}$, we have
\beqn
\int_{B_{\delta R}}g dm  &\leq & \left(\sup\limits_{B_{\delta R}}g \right) m(B_{\delta R})  \leq  \left(\sup\limits_{B_{\delta R}}g^{\frac{1}{2}}\right)^2 m(B_{\delta R})\\
&\leq& e^{2C(1+\sqrt{|K|R})}(1-\delta)^{-2\nu}m(B_{R})^{-1}\left(\int_{B_R}g^{\frac{1}{2}}dm\right)^2 ,
\eeqn
which means that (\ref{alpha}) holds for $\alpha_0=1$, $\alpha=\frac{1}{2}$ and $\sigma=1$. On the other hand, in this case, (\ref{h2}) means that (\ref{Lcondi2}) holds. Then (\ref{U}) becomes
\begin{equation}\label{e-v}
\int_{B_{\delta R}} u dm \leq e^{C_1(1+ \sqrt{|K|}R)} m(B_{R})e^{\bar{v}}.
\end{equation}
Also, letting $g=e^{\bar{v}}u^{-1}$ and taking $\alpha_0=\infty$, $\alpha=1$ and $\sigma=1$,  by the same argument, we can conclude from Lemma \ref{measure} and Lemma \ref{supu-1} that
\begin{equation}\label{e-u}
\sup\limits_{B_{\delta R}} \{u^{-1}\} \leq e^{C_2(1+ \sqrt{|K|}R)}e^{-\bar{v}}.
\end{equation}
Therefore, from $\inf\limits_{B_{\delta R}} u=\left(\sup\limits_{B_{\delta R}} \{u^{-1}\}\right)^{-1}$, we obtain the following
$$
\int_{B_{\delta R}} u dm \leq e^{C_3(1+ \sqrt{|K|}R)} m(B_{R})\inf\limits_{B_{\delta R}}u,
$$
namely, (\ref{mean2-2}) holds. Then we immediately obtain by (\ref{mean2-1})
\begin{equation}
\sup\limits_{B_{\rho \delta R}}u \leq e^{C_4(1+ \sqrt{|K|}R)} \frac{m(B_{R})}{m(B_{\delta R})}\inf\limits_{B_{\delta R}}u\leq e^{C_5(1+ \sqrt{|K|}R)}\inf\limits_{B_{\rho\delta R}}u,
\end{equation}
which is the desired inquality by replacing $\rho\delta$ with $\delta$.
\qed

\vskip 2mm

The following corollary follows directly from Theorem \ref{mean2} by taking $R\rightarrow \infty$.
\begin{cor}{\label{cmean2}}
Let $(M, F, m)$ be an $n$-dimensional forward complete non-compact Finsler measure space with finite reversibility $\Lambda$. Assume that ${\rm Ric}_{\infty} \geq 0$ and $|\tau|\leq k$ for some $k>0$. If $u$ is a positive harmonic function on $M$, then there exist positive constant $C=C\left(n, k, \Lambda\right)$ such that
$$
\sup\limits_{M}u\leq C~\inf\limits_{M}u.
$$
\end{cor}

\vskip 2mm

Finally, we will prove Theorem \ref{ggeh} in the following.
\vskip 2mm
{\it Proof of Theorem \ref{ggeh}.}  Let $u$ is a positive harmonic function, $\Delta u=0$ in a weak sense on $M$. It was proved that $u\in W^{2,2}_{loc}(M)\bigcap \mathcal{C}^{1,\alpha}(M)$(\cite{OHTA}). It follows from the Bochner formula (\ref{BWforinf}) and ${\rm Ric}_{\infty}\geq -K$ $(K\geq 0)$ that
\begin{equation}
\int_M\phi\Delta^{\nabla u} F^2(\nabla u)dm\geq -2 K \int_M \phi F^2(\nabla u) dm
\end{equation}
for all nonnegative functions $\phi \in \mathcal{C}_0^{\infty}\left(M\right)$, which is just (\ref{hF-1}). Hence, by Lemma \ref{Delta} with $p=1$, we have
\be
\sup\limits_{B_{\frac{1}{8}}(x)}F^2(\nabla u)\leq e^{C(1+\sqrt{K})}(1+K)^{\frac{\nu}{2}}m(B_{\frac{1}{4}}(x))^{-1}\int_{B_{\frac{1}{4}}(x)}F^2(\nabla u)dm, \label{g-1}
\ee
where $C=C\left(n, \nu, k, \kappa, \kappa^*\right)>0$ is a universal constant.

In the following, we continue to denote by $C>0$ some universal constant, which may be different line by line. Let $\phi$ a cut-off function defined by
$$
\phi(y)= \begin{cases}1 & \text { on } B_{\frac{1}{4}}(x), \\ 2-4d \left(x, y\right) & \text { on } B_{\frac{1}{2}}(x) \backslash B_{\frac{1}{4}}(x), \\ 0 & \text { on } B_{{1}}(x) \backslash B_{\frac{1}{2}}(x).\end{cases}
$$
Then $F^*(-d \phi) \leq 4$ and $F^*(d \phi) \leq 4\kappa$ a.e. on $B_{\frac{1}{2}}(x)$. It follows that
\beqn
\int_{M}\phi^2 F^2(\nabla u) dm &=&\int_{M}\phi^2 du(\nabla u) dm\\
&=&\int_{M}d(u\phi^2)(\nabla u)dm-2\int_{M}u\phi d\phi(\nabla u)dm\\
&\leq &~\frac{1}{2}\int_{M}\phi^2 F^2(\nabla u) dm +2\int_{M}u^2 F^{*2}(-d\phi) dm,
\eeqn
that is,
$$
\int_{M}\phi^2 F^2(\nabla u) dm\leq 4\int_{M}u^2 F^{*2}(-d\phi) dm.
$$
Therefore,
$$
\int_{B_{\frac{1}{4}}(x)}F^2(\nabla u) dm\leq 64\int_{B_{\frac{1}{2}}(x)}  u^2 dm\leq 64~m\left({B_{\frac{1}{2}}(x)}\right)\left(\sup\limits_{B_{\frac{1}{2}}(x)} u\right)^2.
$$
Combining the above equation with (\ref{g-1}) yields
\begin{equation}\label{g-2}
\sup\limits_{B_{\frac{1}{8}}(x)}F^2(\nabla u)\leq e^{C(1+\sqrt{K})}(1+K)^{\frac{\nu}{2}}\left(\sup\limits_{B_{\frac{1}{2}}(x)}u\right)^2,
\end{equation}
where we have used the volume comparison in Proposition \ref{volcom1}(2). By applying the Theorem \ref{mean2}, we obtain
$$
\begin{aligned}
\sup\limits_{B_{\frac{1}{8}}(x)}F^2(\nabla u) %e^{C(1+\sqrt{|K|})}(\Lambda^2+\mathcal{A})^{\frac{\nu}{2}}\left(\inf\limits_{B_{\frac{1}{2}}(x)}u\right)^2\\
\leq e^{C(1+\sqrt{K})}(1+K)^{\frac{\nu}{2}}\left(\inf\limits_{B_{\frac{1}{8} }(x)}u\right)^2,
\end{aligned}
$$
namely,
$$
F(\nabla u)(x)\leq e^{C(1+\sqrt{K})}(1+K)^{\frac{\nu}{4}}u(x)
$$
for any $x\in M$, where $\nu=4(n+4k)-2$. For $F(\nabla (-\log u))=\overleftarrow{F}(\overleftarrow{\nabla}\log u)$, the same argument works. This means that the gradient estimate (\ref{gradient-1}) holds.
\qed

\end{document}